\documentclass[times]{article}

\usepackage{moreverb}

\usepackage{hyperref}
\usepackage{graphicx}       
\usepackage{longtable} 
\usepackage{amssymb,amsmath}
\usepackage{subfigure}
\usepackage{units}
\usepackage{nicefrac}
\usepackage[utf8]{inputenc}
\usepackage{epsfig}
\usepackage{longtable}
\usepackage{siunitx}
\usepackage{multirow}
\usepackage{units}
\usepackage{nicefrac}
\newcommand{\dd}[2]{\frac{d}{d#1}{\left( #2 \right)}}
\newcommand{\ns}{\mathcal H_K(\Omega)}
\newcommand{\R}{\mathbb R}

\begin{document}

\author{Tobias K\"oppl$^1$ \corrauth, Gabriele Santin$^2$, Bernard Haasdonk$^2$, Rainer Helmig$^1$}

\title{Numerical modelling of a peripheral arterial stenosis using dimensionally reduced models and kernel methods}
 
\author{Tobias K\"oppl, Rainer Helmig \\ \\
 { \small University of Stuttgart} \\
 {\small Department of Hydromechanics and Modelling of Hydrosystems} \\
 {\small Pfaffenwaldring 61} \\
 \\ \\
 Gabriele Santin, Bernard Haasdonk \\ \\
 { \small University of Stuttgart} \\
 {\small Institute of Applied Analysis and Numerical Simulation} \\
 {\small Pfaffenwaldring 57}  
}

% \author{T. K\"oppl  \and 
% 	G. Santin	\and
% 	B. Haasdonk \and
% 	R. Helmig
% 	}
% 	
\maketitle
% \address{\begin{center}
%           $^1$University of Stuttgart, Department of Hydromechanics and Modelling of Hydrosystems, Pfaffenwaldring 61, \\ 
%          D-70569 Stuttgart, Germany \\         
%          $^2$University of Stuttgart, Institute of Applied Analysis and Numerical Simulation, Pfaffenwaldring 57,\\
%          D-70569 Stuttgart Germany 
%          \end{center}}

%\corraddr{tobias.koeppl@iws.uni-stuttgart.de}

\begin{abstract}
In this work, we consider two kinds of model reduction techniques to simulate blood flow through the largest systemic arteries, where a stenosis is located
in a peripheral artery i.e. in an artery that is located far away from the heart. For our simulations we place the stenosis in one of the tibial arteries 
belonging to the right lower leg (right post tibial artery). The model reduction techniques that are used are on the one hand dimensionally reduced models
(1-D and 0-D models, the so-called mixed-dimension model) and on the other hand surrogate models produced by kernel methods. Both methods are combined in 
such a way that the mixed-dimension models 
yield training data for the surrogate model, where the surrogate model is parametrised by the degree of narrowing of the peripheral stenosis. By means of a 
well-trained
surrogate model, we show that simulation data can be reproduced with a satisfactory accuracy and that parameter optimisation or state estimation problems can 
be solved in a very efficient way.
Furthermore it is demonstrated that a surrogate model enables us to present after a very short simulation time the impact of a varying degree of stenosis on 
blood flow, obtaining a speedup of several orders over the full model.
\end{abstract}

% \keywords{blood flow simulations, peripheral stenosis, dimensionally reduced models, mixed-dimension models, kernel methods, surrogate models,
% real-time simulations}

%\maketitle

\section{Introduction}

During the recent decades, the interest in numerical simulation of blood flow has been growing continuously. Main reasons for this development are
the increase of computational power, the design of efficient numerical algorithms and the improvement of imaging techniques combined
with elaborated reconstruction techniques yielding data on the geometry of interesting objects as well as important modelling parameters like
densities of a fluid or tissue \cite{ambrosi2012modeling,formaggia2010cardiovascular,quarteroni2000computational}. 
The motivation for putting more and more effort into these developments has been evoked 
by the fact that computational tools enable clinical doctors and physiologists 
to obtain some insight into cardiovascular diseases in a non-invasive way. By this, the risk of infections and other dangers
can be remarkably reduced. Using numerical models, it is, e.g., possible to make predictions on how a stenosis affects the blood supply of organs, in 
particular
one wants to find out to what extent a vessel can be occluded without reducing the blood flow significantly. 
In this context, it is also of great interest whether blood vessel systems like the 
Circle of Willis \cite{alastruey2007modelling,ryu2015coupled} or small interarterial connections \cite{schaper1999collateral,scholz2002contribution}
can help to restore the reduced blood flow. Furthermore, the simulation of fundamental regulation mechanisms like vasodilation,
arteriogenesis and angiogenesis help to understand how the impact of an occlusion on blood flow can be reduced
\cite{alastruey2008reduced,drzisga2016numerical,koppl2014influence}. Further important application areas are the 
stability analysis of an implanted stent or an aneursysm. Thereby, it is crucial to compute realistic pressure values, mass fluxes and wall shear 
stresses in an efficient way such that these data can be evaluated as fast as possible. 

However, the simulation of flow through a cardiovascular system
is a very complex matter. Since it is composed of a huge number of different vessels that are connected in a very complex way, it is usually not possible to 
resolve large parts of the cardiovascular systems due to an enormous demand for computational power and data volume. In addition to that there are many 
different kinds of vessels covering a large range of radii, wall thicknesses and lengths \cite[Chap. 1, Tab. 1.1, Tab. 1.2]{formaggia2010cardiovascular}. 
The vessel walls of arteries are, e.g., much thicker than those of veins, due to the fact that they have to transport blood at high pressure. 
As the heart acts like a periodic pump, blood flow in the larger arteries and venes having an elastic vessel wall is pulsatile and exhibits high Reynolds numbers 
\cite[Chap. 1, Tab. 1.7]{formaggia2010cardiovascular}. This requires the usage of FSI algorithms and discretisation techniques 
for convection dominated flows \cite{bazilevs2009patient,colciago2014comparisons,crosetto2011fluid,torii2009fluid}. 
Contrary to that flows in small vessels with respect to diameter or length \cite[Chap. 1, Tab. 1.7]{formaggia2010cardiovascular} 
exhibit small Reynolds numbers. Moreover, the walls of such vessels are not significantly deformed and therefore they
can be modelled as quasi-rigid tubes. Vessels having these properties can be typically classified as arterioles, venoles or capillaries.
Taking all these facts into account, it becomes obvious that to this part of the cardiovascular system totally different models and methods have to be applied
\cite[Chap. 6.2]{d2007multiscale}\cite{notaro2016mixed,reichold2009vascular}.

Due to the variety of different vessel geometries and types of flows, it is unavoidable to consider a coupling of different kinds of models such that realistic
blood flow simulations within the entire or within a part of the cardiovascular system can be performed. In order to reduce the complexity of the numerical model
mixed-dimension models have been introduced \cite{alastruey2008lumped,alastruey2012physical,formaggia2006numerical,malossi2013implicit,marchandise2009numerical,
marsden2015multiscale,mynard20081d,passerini20093d}. 
Thereby, subnetworks of larger vessels are modelled by three-dimensional (3-D) or one-dimensional (1-D) PDEs in space. At the inlets and outlets of these networks, 
the corresponding models are coupled with one-dimensional (1-D) PDEs or zero-dimensional (0-D) models (systems of ODEs) incorporating e.g. the windkessel effect
of the omitted vessels and the pumping of the heart \cite{vcanic2003mathematical,ismail2013adjoint,olufsen1999structured}. 
An alternative to these open-loop models for arterial or venous subnetworks of the systemic circulation 
are closed-loop models linking the inlets and the outlets of a subnetwork by a sequence of 
0-D models for the organs, pulmonary circulation and the heart \cite{liang2005closed,liang2009multi,mynard2015one}. 

Besides the usage of dimensionally reduced or
multiscale models a further method of reducing the complexity of blood flow simulations has been established in the recent years
\cite{manzoni2012model,quarteroni2003optimal,quarteroni2014reduced}. 
This approach is called Reduced Basis method (RB method) and is based on the idea to solve in an offline phase
parameterised PDEs for a few parameters (see e.g. \cite{Haasdonk2017}). These so-called solution snapshots are then used as basis of a low-dimensional space, 
and projecting the PDE to this 
space results in a low-dimensional problem. This reduced system can be solved in a so-called online phase, where solutions for multiple parameters can be 
efficiently computed. Within relevant biomedical application areas these parameters determine usually the
shape of a bifurcation, a stenosis, a bypass or an inflow profile \cite[Chap. 8]{quarteroni2015reduced}, \cite{quarteroni2003optimal,quarteroni2014reduced}.

In this paper, we want to investigate the performance of a different type of surrogate model obtained via machine learning techniques, and in particular with kernel 
methods 
\cite{Fasshauer2015,haasdonk2005transformation,haasdonk2007invariant,smola1998learning,Wendland2005}. 
Contrary to the RB method, the surrogate model is in this case represented by a linear combination of kernel functions like the Gaussian or the Wendland kernel 
\cite{Wendland1995a,wirtz2015surrogate}. The coefficients in the linear combination and the parameters for the kernel functions are obtained from a 
training and validation process which is performed in an offline or training phase \cite{HS2017a}. These methods have the advantage of constructing 
nonlinear, data-dependent surrogates that can reach significant degrees of accuracy, while not needing an excessively large amount of data, as is instead the case for 
different machine learning techniques.

Using this method we want to simulate the impact of a stenosis on blood flow, in particular we consider a peripheral arterial stenosis. This type of
stenosis is of high interest for physiologists, since peripheral arteries supply organs that are located far away from the heart. It is obvious that peripheral organs
and tissues are affected by a potential risk of undersupply and therefore an occlusion of peripheral arteries is extremely critical. For our simulations, we place 
a stenosis in the right tibial artery located in the lower right leg (see Figure \ref{fig:ArterialNetwork_Stenosis})
and study the pressure and flow rate curves, i.e., the evolution of the two quantities at different points over a complete heart-beat, in the vicinity of this 
stenosis. 

We use for the blood flow 
simulations a 1D-0D coupled model, where we assign a 1-D flow model to the $55$ main arteries of the systemic circulation 
\cite{alastruey2008reduced,alastruey2008lumped,d2007multiscale,koppl2013reduced}. At the outlets of this network 0-D models are attached to the 1-D models to incorporate
the influence of the omitted vessels. The stenosis is modelled by a 0-D model consisting of an ODE
\cite{mynard20081d,seeley1976effect,young1973flow}. Using the dimensionally reduced
model, we can produce realistic pressure and flow rate curves in a fast way. 

Although this model is already both accurate and relatively fast, it is still too slow and hence not suitable for real-time simulation or parameter estimation. 
To overcome this problem, we train a kernel-based surrogate model which predicts, depending on the degree of stenosis, a pressure or flow-rate curve.
The surrogate is trained in a data-dependent way by computing pressure and flow rate curves for different
degrees of stenosis, which are used as training data for the kernel method to construct an accurate
surrogate model. The intention of this modelling approach is that a
combination of dimensionally reduced models and kernel methods allows us to simulate the impact of a stenosis for an
arbitrary degree of narrowing in a very short time. 
Simulation techniques of this kind might support clinical doctors and researchers with some important information after a relative short time,
such that their diagnostic process can be optimised.

This combination of techniques is relatively new, and the results presented in this paper demonstrate its effectiveness. Moreover, although we concentrate 
here on the prediction of pressure and flow rate curves, the same technique can be easily adapted to construct surrogates of other relevant quantities of the blood flow 
simulation.

We remark that the present approach has potentially different advantages over the RB method. Indeed, the RB method typically requires the 
computation of 
several time snapshots in the offline phase to simulate the time evolution, and also requires a time integration, with mostly the same timestep used in the 
full model, during 
the online phase. On the other hand, kernel based surrogates only require the time evolution of the quantities of interest in the desired time 
interval as training data, and in the online phase they can directly predict it for a new parameter, without the need of any time integrations. Moreover, 
it is well known that the RB 
method may perform poorly when applied to transport problems, especially in the presence of moving structures or discontinuities evolution. This behavior is reflected in 
a slowly decreasing Kolmogorov width, as it is discussed e.g. in \cite[Example 3.4]{Haasdonk2017}.

The remainder of our work is structured as follows: In Section \ref{sec:DimRedModel} we outline all the details of the 1D-0D coupled model and the model for the 
stenosis. In addition to that some comments on the numerical methods are presented. The first subsection of the following Section \ref{sec:KerMeth} contains some 
information on the fundamentals of kernel methods. This first part presents results that are mainly already discussed in the cited literature. Nevertheless, 
since we aim to address researcher of both the blood flow and machine learning communities, we include it to provide a clearer explanation for the interested reader, 
which may be acquainted on only one of the two fields.
The second and third subsection describe how kernel methods can be used to compute flow variables in 
dependence
of the degree of stenosis. By means of the models and methods from Section \ref{sec:DimRedModel} and \ref{sec:KerMeth}, we perform in Section \ref{sec:NumTests} some
numerical tests illustrating the accuracy of the surrogate kernel model. Moreover it is shown how the surrogate kernel model can be used to solve a 
state estimation problem. The paper is concluded by Section \ref{sec:Concl}, in which we summarise the main results and make some comments on possible future 
work.

\section{Simulation of arterial blood flow by dimensionally reduced models}
\label{sec:DimRedModel}

\begin{figure}[h!]
\begin{center}
\includegraphics[scale=0.35]{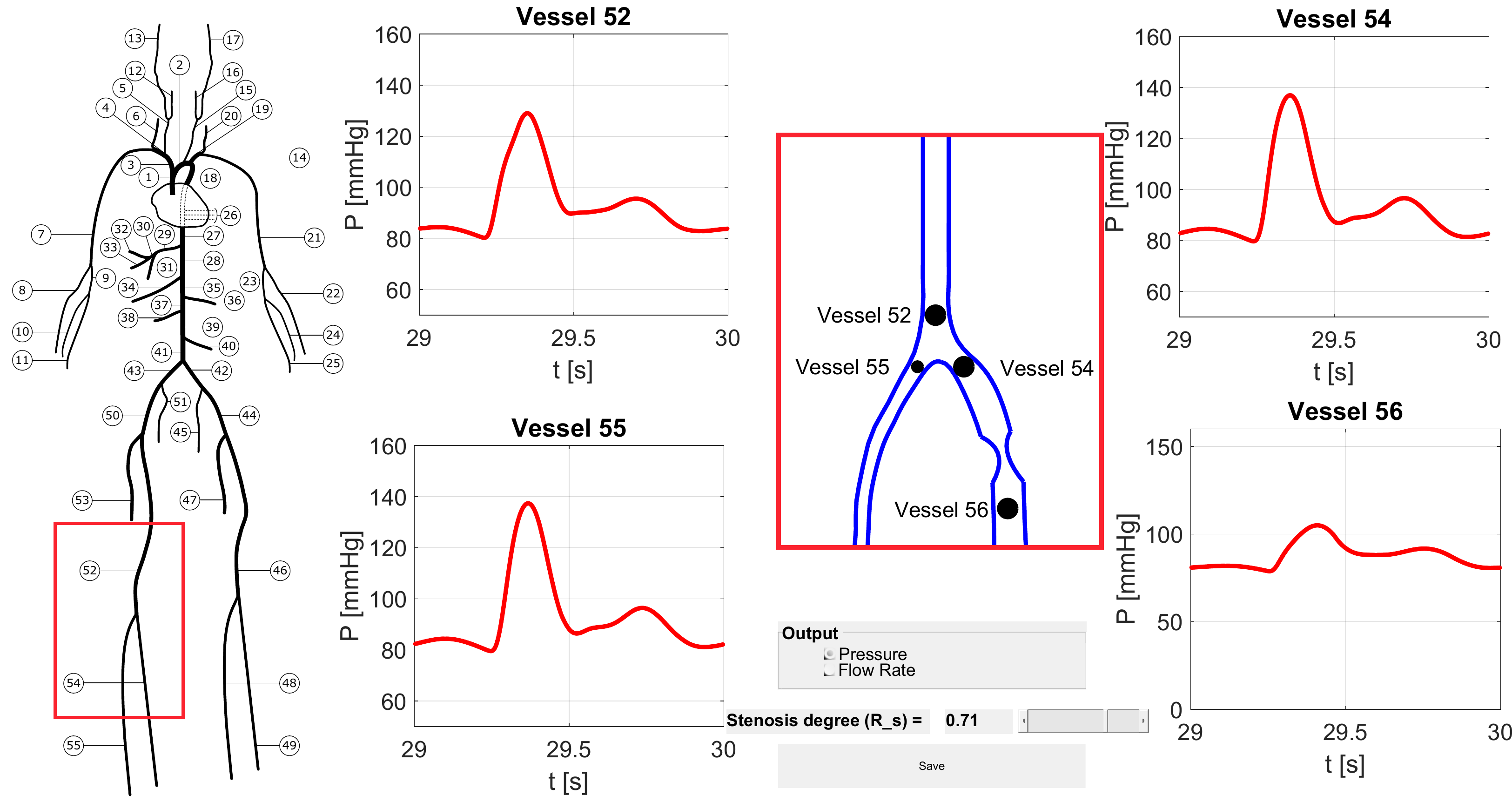}
\caption{\label{fig:ArterialNetwork_Stenosis} This figure shows an arterial network consisting of the $55$ main arteries of the systemic circulation
\cite{sherwin2003computational,wang2004wave,westerhof1969analog}. In Vessel $54$ (right anterior tibial artery) we put a stenosis and study the effect on blood flow. The 
stenosis
splits Vessel $54$ into two parts. The proximal part is again labelled with the index $54$, while the distal part receives the index $56$. At the places of the
black dots that are located at the outlet of Vessel $52$, at the inlets of Vessel $54$ and Vessel $55$ as well as at the inlet of 
Vessel $56$, we report over one heart beat pressure and flow rate curves. Samples of these curves serve as training data for the kernel methods creating 
a surrogate model which maps the degree of stenosis $R_s \in \left[0,1 \right]$ to the samples of the corresponding curves.}
\end{center}
\end{figure}

Simulating
blood flow from the heart to the arms and legs, we consider the arterial network, presented in \cite{sherwin2003computational,wang2004wave,westerhof1969analog}. This network consists
of the $55$ main arteries including the aorta, carotid arteries, subclavian arteries and the tibial arteries (see Figure \ref{fig:ArterialNetwork_Stenosis}). 
Our modelling approach for simulating blood flow through this network is based on the idea to decompose in a first step 
the network into its single vessels. In a next step a simplified 1-D flow model is assigned to each vessel. Finally, 
the single models have to be coupled at the different interfaces, in order to obtain global solutions for the flow variables.

The following subsections present the basic principles of the 1-D model and the coupling conditions
at bifurcations as well as at the stenosis. Furthermore, we make some comments on the numerical methods that are used to compute a suitable solution.
At the inlet of the aorta (Vessel $1$, Figure \ref{fig:ArterialNetwork_Stenosis} left), we try to emulate the 
heart beats by a suitable boundary condition. In order to account for the windkessel effect of the omitted vessels, the 1-D flow models associated with the 
terminal vessels are coupled with ODE-systems (0-D models), which are derived from electrical science \cite{alastruey2008lumped}. 
Usually the term windkessel effect is related to the ability of large deformable vessels
to store a certain amount of blood volume such that a continuous supply of organs and tissue can be ensured. However, also many of the vessels that are not 
depicted in Figure
\ref{fig:ArterialNetwork_Stenosis} exhibit this feature to some extent. Furthermore the arterioles located beyond the outlets of this network can impose some resistance
on blood flow \cite{alastruey2008reduced,koppl2014influence}. These features have to be integrated into the outflow models in order to be able to 
simulate realistic pressure and flow rate curves. 

\subsection{Modelling of blood flow through a single vessel}
\label{sec:SingVess}

Let us suppose that the Navier-Stokes equations are defined on cylindrical domains of length $l_i\in\mathbb{R},\;i \in \left\{1,\ldots,55 \right\}$ and that 
their main axis 
are aligned with
the $z$-coordinate. Modelling the viscous forces, we assume that blood in large and medium sized arteries can be treated as an incompressible
Newtonian fluid, since blood viscosity is almost constant within large and middle sized vessels \cite{pries1992blood}. The boundaries of 
the computational domain change in time due to the elasticity of the arterial vessels walls and the pulsatile flow. 
Thereby, it is assumed that the vessel displaces only in the 
radial direction and that the flow is symmetric with respect to the main axis of the vessel \cite{formaggia2006numerical,quarteroni2004mathematical}. In addition to
that we postulate that the $z$-component $\mathbf{u}_z$ of the velocity flield $\mathbf{u}$ is dominant with respect to the other components. Taking all these 
assumptions into account and integrating the Navier-Stokes equations across a section area $S\left(z,t \right)$ perpendicular to the $z$-axis at the 
place $z \in \left(0,l_i \right)$ and a time $t>0$, one obtains the following system of PDEs \cite{barnard1966theory,vcanic2003mathematical,hughes1973one}:
\begin{align}
\label{eq:masscon} 
\frac{\partial A_i}{\partial t} + \frac{\partial Q_i}{\partial z} = 0,&\; z \in \left(0,l_i \right),\;t>0, \\
\nonumber
\\
\label{eq:momentumcon}
\frac{\partial Q_i}{\partial t} + \frac{\partial}{\partial z}\left( \frac{Q_i^2}{A_i} \right) + \frac{A_i}{\rho}\frac{\partial p_i}{\partial z} = -K_{r} \frac{Q_i}{A_i},
&\; z \in \left(0,l_i \right),\;t>0.
\end{align}
The unknowns of this system $A_i$, $Q_i$ and $p_i$ 
denote the section area of vessel $V_i$,
the flow rate and the averaged pressure within this vessel. Mathematically, these quantities are defined by the following integrals:
$$
A_i \left(z,t \right) = \int_{S\left(z,t \right)} \;dS, \quad Q_i \left(z,t \right) = \int_{S\left(z,t \right)} \mathbf{u}_z\;dS, \quad
p_i\left(z,t \right) = \frac{1}{A_i} \int_{S\left(z,t \right)} P_i\;dS,
$$
where $P_i$ is the 3-D pressure field. Please note that \eqref{eq:masscon} is the 1-D version of the mass conservation equation, 
while \eqref{eq:momentumcon} represents the 1-D version of the momentum equations. $\rho$ stands for the density of blood. 
$K_{r}$ is a resistance parameter containing the dynamic viscosity $\eta$ of blood 
\cite{smith2002anatomically}:
$K_r = 22 \pi \eta/\rho$. The PDE system \eqref{eq:masscon} and \eqref{eq:momentumcon} is closed by means of an algebraic equation, which can be derived from the 
Young-Laplace 
equation \cite{olufsen1999structured}:
\begin{equation}
\label{eq:PressureAreaRelation}
p_i(z,t) = G_{0,i} \left( \sqrt{\frac{A_{i}}{A_{0,i}}} - 1 \right),\quad G_{0,i} = \frac{\sqrt{\pi} \cdot h_{0,i} \cdot
E_i}{\left(1-\nu^2\right) \cdot \sqrt{A_{0,i}}},
\end{equation}
where $E_i$ is the Young modulus, $A_{0,i}$ stands for the section area at rest,
$h_{0,i}$ is the vessel thickness and $\nu$ is the Poisson ratio. 
Due to the fact that biological tissue is practically incompressible, $\nu$ is
chosen as follows: $\nu = 0.5$. Equation \eqref{eq:PressureAreaRelation} assumes that the vessel wall is instantaneously in equilibrium with the forces acting on it. 
Effects like wall inertia and viscoelasticity could be incorporated using a differential pressure law \cite{devault2008blood,mynard2015one,valdez2009analysis}.
However, neglecting the viscoelasticity, maintains the strict hyperbolicity of the above PDE system \cite{Riviere}. Therefore,
the interaction between the blood flow and the elastic vessel walls \eqref{eq:PressureAreaRelation} is accounted for by using \eqref{eq:PressureAreaRelation}.
Assuming that $G_{0,i}$ and $A_{0,i}$ are constant, the PDE-system \eqref{eq:masscon}-\eqref{eq:PressureAreaRelation} can be represented in a compact form:
\begin{equation}
\label{eq:PDESys}
\frac{\partial \mathbf{U}_i}{\partial t} + \frac{\partial \mathbf{F}_i}{\partial z}\left( \mathbf{U}_i \right) = \mathbf{S}_i\left( \mathbf{U}_i \right),
\;z \in \left(0,l_i \right),\;t>0.
\end{equation}
For $\mathbf{U}_i = \left( A_i, Q_i \right)^T$, the flux function $\mathbf{F}_i$ and the source function $\mathbf{S}_i$ are given by:
$$
\mathbf{F}_i \left( \mathbf{U}_i \right) = \begin{pmatrix}
                                            Q_i \\ \frac{Q_i^2}{A_i} + \frac{A_i^{\frac 32}}{\rho \sqrt{A_{i,0}}}
                                            \end{pmatrix} \;\text{ and }\;
\mathbf{S}_i \left( \mathbf{U}_i \right) = \begin{pmatrix}   
                                            0 \\ - K_{r,i} \frac{Q_i}{A_i}
                                           \end{pmatrix}.                                    
$$
This system may be written in a quasilinear form:
\begin{equation*}
\frac{\partial \mathbf{U}_i}{\partial t} + \nabla_{\mathbf{U}_i}\mathbf{F}_i\frac{\partial \mathbf{U}_i}{\partial z}
= \mathbf{S}_i\left( \mathbf{U}_i \right),\;z \in \left(0,l_i \right),\;t>0,
\end{equation*}
where $\nabla_{\mathbf{U}_i}\mathbf{F}_i$ is the $2 \times 2$ Jacobian matrix of the flux function $\mathbf{F}_i$, having the eigenvalues
$\lambda_{i,1}$ and $\lambda_{i,2}$. Denoting by $v_i = Q_i/A_i$ the fluid velocity and $v_{c,i} \left(A_i \right)$ the characteristic wave velocity of vessel $V_i$,
it can be shown that: $\lambda_{i,1} = v_i -  v_{c,i}$ and $\lambda_{i,2} = v_i + v_{c,i}$. Under physiological conditions, it can be observed that 
\cite{formaggia2002one}:
\begin{equation}
\label{eq:wavespeed}
v_i = \frac{Q_i}{A_i} \ll \sqrt{\frac{G_{0,i}}{2 \rho} \sqrt{\frac{A_i}{A_{0,i}}}} = v_{c,i} \left(A_i \right).
\end{equation}
Therefore it holds for the eigenvalues: $\lambda_{i,1} < 0$ and $\lambda_{i,2}>0$ and the above PDE-system is hyperbolic. Exploiting the fact that
the Jacobian matrix $\nabla_{\mathbf{U}_i}\mathbf{F}_i$ is diagonalisable, there is an invertible matrix $L_i$ such that it can be decomposed as follows:
$\nabla_{\mathbf{U}_i}\mathbf{F}_i = L_i^{-1} \Lambda_i L_i$, where $\Lambda_i$ is a diagonal matrix that has the eigenvalues $\lambda_{i,1}$ and $\lambda_{i,2}$ on
its diagonal. By this, the PDE-system \eqref{eq:PDESys} can written in its characteristic variables
$\mathbf{W}_i = \left( W_{1,i}, W_{2,i} \right)^T$:
\begin{equation}
\label{eq:charSys}
\frac{\partial \mathbf{W}_i}{\partial t} + \Lambda_i \frac{\partial \mathbf{W}_i}{\partial z} = L_i \mathbf{S}_i\left( \mathbf{W}_i \right),
\;z \in \left(0,l_i \right),\;t>0.
\end{equation}
The characteristic variables $\mathbf{W}_i$ and $L_i$ are related by the following equation: 
\begin{equation}
\label{eq:ODEW}
\frac{ \partial \mathbf{W}_i }{ \partial \mathbf{U}_i } = L_i, \; \mathbf{W}_i(\mathbf{U}_i) = \begin{pmatrix} 0 \\ 0  \end{pmatrix}\; \text{ for } \;
\mathbf{U}_i = \begin{pmatrix} A_{0,i} \\ 0  \end{pmatrix}.
\end{equation}
An integration of \eqref{eq:ODEW} 
yields that the characteristic variables  $W_{1,i}$ and $W_{2,i}$ can be expressed by the primary variables $A_i$ and $Q_i$ as follows: 

\begin{align}
\label{eq:W1} 
W_{1,i} &= -\frac{Q_i}{A_i} + 4 \sqrt{ \frac{G_{0,i}}{2 \rho} } \left( \left( \frac{A_i}{ A_{0,i} } \right)^{\frac 14} - 1 \right) = -v_i + 4 \cdot \left( v_{c,i} 
\left(A_i \right)
- v_{c,i}\left(A_{0,i} \right) \right), \\  
\nonumber
\\
\label{eq:W2}
W_{2,i} &= \frac{Q_i}{A_i}  + 4 \sqrt{ \frac{G_{0,i}}{2 \rho} } \left( \left( \frac{A_i}{ A_{0,i} } \right)^{\frac 14} - 1 \right) = v_i + 4 \cdot \left( v_{c,i} 
\left(A_i \right)
- v_{c,i}\left(A_{0,i} \right) \right).
\end{align} 
Based on condition \eqref{eq:wavespeed} and the signs of $\lambda_{1,i}$ and $\lambda_{2,i}$, 
one can prove that $W_{1,i}$ is a backward and $W_{2,i}$ is a forward travelling wave \cite{vcanic2003mathematical}
\cite[Chap. 2]{d2007multiscale}. Furthermore, it can be shown that $W_{1,i}$ and $W_{2,i}$ are moving on characteristic curves $c_{j,i}$
defined by two ODEs:
\begin{equation}
\label{eq:CharCurves}
\frac{d c_{j,i}}{dt} \left(t \right) = \lambda_j\left( c_{j,i}(t),t \right),\;j \in \left\{1,2 \right\}.
\end{equation}
These insights are crucial for a consistent coupling of the submodels at the different interfaces, since it reveals that at each inlet and
outlet of a vessel exactly one coupling or boundary condition has to be imposed. The other condition is obtained from the outgoing characteristic variable. At the
inlet $z=0$, the variable $W_{1,i}$ is leaving the computational domain, whereas at $z=l_i$ the variable $W_{2,i}$ is the outgoing characteristic variable.

\subsection{Numerical solution techniques}
\label{sec:NMC}

According to standard literature \cite[Chap. 13]{quarteroni2010numerical}, the main difficulties that arise in terms of numerical treatment of hyperbolic PDEs are to 
minimise
dissipation and dispersion errors, in order to avoid an excessive loss of mass and a phase shift for the travelling waves. A standard remedy for these problems is to 
apply
higher order discretisation methods in both space and time \cite{sherwin2003computational} such that the numerical solution is as accurate as possible. However, higher 
order discretisation methods tend to create oscillations in the vicinity of steep gradients or sharp corners, which can be removed by some additional postprocessing
\cite{gottlieb2001strong,koppl2013reduced,krivodonova2007limiters}. 
Moreover, time stepping methods of higher order require small time steps to be able to
resolve the dynamics of a fast and convection dominated flow and to fulfill a CFL-condition, if they are explicit. 

Considering all these features, we use in this work the numerical method of characteristics (NMC), which is 
explicit and of low approximation order (first order in space and time \cite[Theorem 1]{Riviere}), 
leading to large dissipation and dispersion errors. This drawback can be circumvented by using a fine grid in space and sufficiently small timesteps.
Since we deal in this work with 1-D problems, a fine grid in space is affordable with respect to computational effort. On the other hand, a fine grid might force an 
explicit method
to exert very small time steps. However, for the NMC it can be proven that its time stepsize is not restricted by a condition of CFL type \cite[Prop. 2]{Riviere}. This 
means
that the NMC can use a fine grid in space and time stepsizes that are small enough to capture the convection dominated blood flow and large enough to 
have an acceptable number of timesteps.

\begin{figure}[h!]
\begin{center}
\includegraphics[scale=0.35]{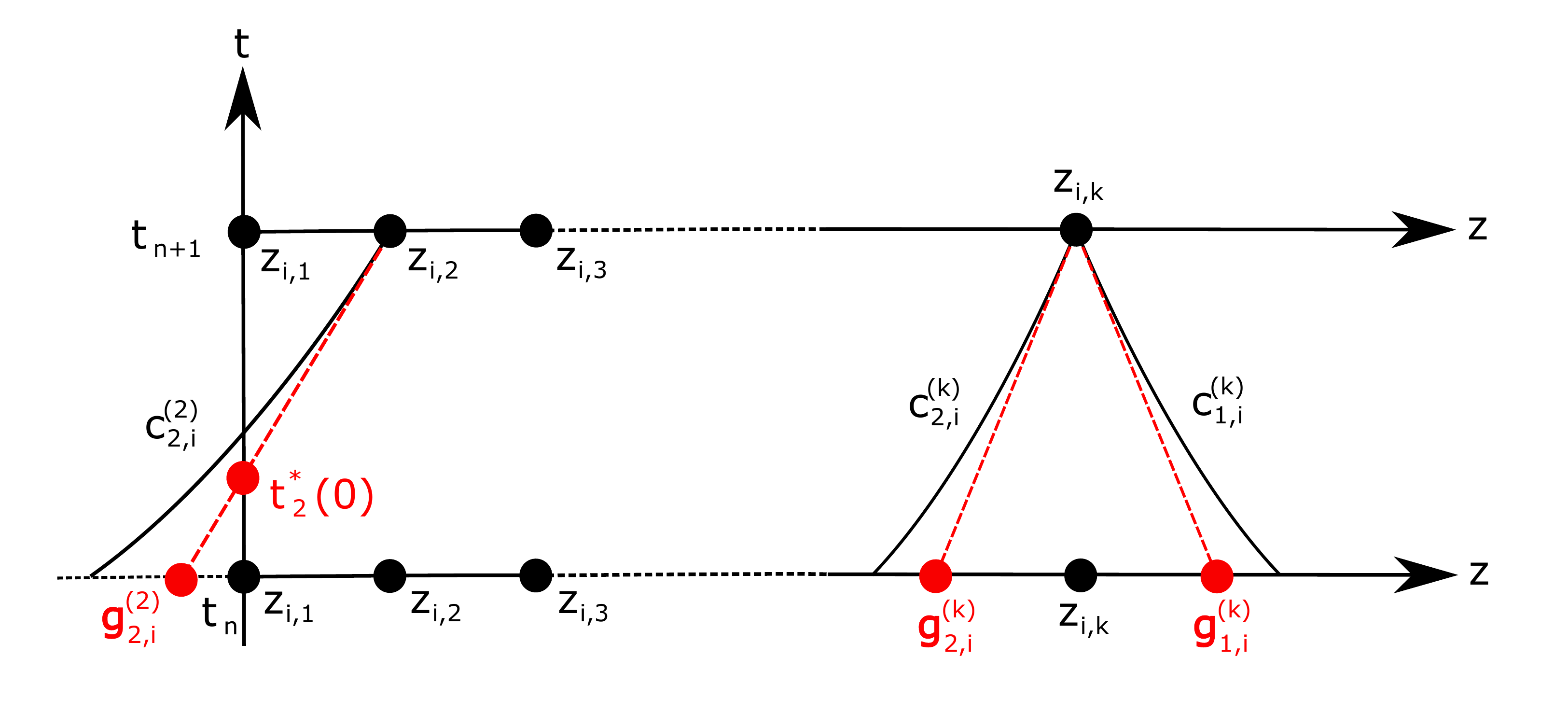}
\end{center}
\caption{\label{fig:NMC_Extrapolation} Linearisation of the characteristic curves $c_{j,i}^{(k)}$ for grid nodes $z_{i,k}$ in the vicinity of
the boundary and in the inner of the computational domain.}
\end{figure}

Let us suppose now that the interval $\left[0,l_i \right]$ for Vessel $V_i$ is discretised by a grid having a meshsize $\Delta z_i$ and grid nodes $z_{i,k} = k 
\cdot \Delta z_i
\in \left[0,l_i \right],\; k \in \left\{0,\ldots,N_{h,i} \right\}$. Here, $N_{h,i}$ is the index of the last grid node. In a time step 
$\left[t_n,t_n+\Delta t \right]$, 
the NMC iterates over all the grid nodes. 
At each grid node $z_{i,k}$ there are two characteristic curves $c_{1,i}^{(k)}$ and $c_{2,i}^{(k)}$ for $W_{1,i}$ and $W_{2,i}$ (see Figure \ref{fig:NMC_Extrapolation}).
Both curves are linearized in $z_{i,k}$ and $t_{n+1} = t_n + \Delta t$. In a next step, the resulting tangents are traced back to the previous time point $t_n$, where 
the corresponding intersection points are denoted by $g_{1,i}^{(k)}$ and $g_{2,i}^{(k)}$, respectively (see Figure \ref{fig:NMC_Extrapolation}). 

This procedure is equivalent to solving for every $z_{i,k}$ the final value problem \cite{Riviere}[Equation (22)], which can be derived from \eqref{eq:CharCurves}:
$$
\frac{d c_{j,i} \left(z_i, t_{n+1}, t \right)}{dt} \ = \lambda_j\left( c_{j,i} \left(z_i, t_{n+1}, t \right), t \right),\; 
c_{j,i} \left(z_i, t_{n+1}, t_{n+1} \right) = z_{i,k}, \;j \in \left\{1,2 \right\}.
$$
Setting $c_{j,i} \left(z_i, t_{n+1}, t_n \right) = g_{j,i}^{(k)}$ we have by a first order approximation:
\begin{align*}
c_{j,i} \left(z_i, t_{n+1}, t_{n+1} \right) - c_{j,i} \left(z_i, t_{n+1}, t_{n} \right)  &= \int_{t_n}^{t_{n+1}} \lambda_j \left( c_{j,i} \left(z_i, t_{n+1}, t \right), t 
\right)\;dt, \\
z_{i,k}-g_{j,i}^{(k)}\left( t_{n} \right) &= \int_{t_n}^{t_{n+1}} \lambda_j \left( c_{j,i} \left(z_i, t_{n+1}, t \right), t \right)\;dt \\
&\approx \Delta t \lambda_j \left( c_{j,i} \left(z_i, t_{n+1}, t_n \right), t_n \right) \\
&= \Delta t \lambda_j \left(  g_{j,i}^{(k)}\left( t_{n} \right), t_n \right). 
\end{align*}
Restricting the PDE-system \eqref{eq:charSys} to the characteristic curves $c_{j,i}^{(k)}$, we have to solve the following ODEs in order to determine approximations for 
$W_{1,i}$ and $W_{2,i}$:
\begin{equation}
\label{eq:ODEWChar}
\frac{d}{dt} \mathbf{W}_{j,i}\left( c_{j,i}^{(k)}(t),t \right) = L_i \mathbf{S}_i \left( \mathbf{W}_{j,i}\left( c_{j,i}^{(k)}(t),t \right) \right).
\end{equation}
An explicit first order discretisation of \eqref{eq:ODEWChar} yields the following extrapolation formula for $\mathbf{W}_{j,i}$ at $c_{j,i}^{(k)}\left(t_{n+1}\right) = 
z_{i,k}$:
\begin{equation}
\label{eq:Extrapolation}
\mathbf{W}_{j,i}\left( z_{i,k},t_{n+1} \right)= \mathbf{W}_{j,i}\left( g_{j,i}^{(k)}\left( t_{n} \right),t_{n} \right) 
+ \Delta t \cdot L_i \mathbf{S}_i \left( \mathbf{W}_{j,i}\left( g_{j,i}^{(k)}\left( t_{n} \right),t_{n} \right) \right).
\end{equation}
At the old time step $t_n$, the values $\mathbf{W}_{j,i}\left( g_{j,i}^{(k)}\left( t_{n} \right),t_{n} \right)$ can be interpolated using the precomputed values at the 
grid nodes
$z_{i,k}$. For large time steps, it may happen that $g_{j,i}^{(k)}\left( t_{n} \right) \notin \left[0,l_i \right]$ (see Figure \ref{fig:NMC_Extrapolation}). 
In these cases, the values 
$\mathbf{W}_{j,i}\left( g_{j,i}^{(k)}\left( t_{n} \right),t_{n} \right)$ can not be interpolated from the spatial values. Therefore, we require a temporal interpolation
at a time point at which the linearised characteristic curves leave the computational domain. At the boundaries $z_{i}=0$ and $z_{i}=l_i$ these time points 
$\left\{ t^\ast_k\left(0 \right),t^\ast_k\left( l_i \right) \right\} \subset \left[t_n,t_{n+1}\right]$ can be computed by:
\begin{equation*}
t^\ast_k\left( 0 \right) = t_n - \Delta t \frac{g_{2,i}^{(k)}}{z_{i,k}-g_{2,i}^{(k)}} \;\text{ and }\;
t^\ast_k\left( l_i \right) = t_n + \Delta t \frac{l_i-g_{1,i}^{(k)}}{z_{i,k}-g_{1,i}^{(k)}}.
\end{equation*}
After computing the ingoing characteristic variables $\mathbf{W}_{2,i}\left( 0, t_{n+1} \right)$ and $\mathbf{W}_{1,i}\left( l_i, t_{n+1} \right)$ for the new time step 
by means of 
an external model (see Subsection \ref{sec:Heart}-\ref{sec:Stenosis}), we use a linear interpolation to provide a surrogate value for the missing characteristic variable:
\begin{align}
\label{eq:WBoundary}
\nonumber
\mathbf{W}_{2,i}\left( 0, t^\ast_k\left( 0 \right) \right) &= \frac{t_{n+1} - t^\ast_k\left( 0 \right)}{\Delta t}\mathbf{W}_{2,i}\left( 0, t_{n} \right) 
+ \frac{t^\ast_k\left( 0 \right)-t_n}{\Delta t} \mathbf{W}_{2,i}\left( 0, t_{n+1} \right),\\
\mathbf{W}_{1,i}\left( l_i, t^\ast_k\left( l_i \right) \right) &= \frac{t_{n+1} - t^\ast_k\left( l_i \right)}{\Delta t}\mathbf{W}_{1,i}\left( l_i, t_{n} \right) 
+ \frac{t^\ast_k\left( l_i \right)-t_n}{\Delta t} \mathbf{W}_{1,i}\left( l_i, t_{n+1} \right).
\end{align}

\subsection{Modelling of heart beats}
\label{sec:Heart}

At the inlet of the aorta, we couple the corresponding 1-D model with a lumped parameter model (0-D model) for the left ventricle of the heart. By means of this
model and the outgoing characteristic variable $W_{1,1}$ the missing ingoing variable $W_{2,1}$ can be determined. In order to compute the pressure $P_{LV}$ in the left 
ventricle, we
consider the following elastance model \cite{formaggia2006numerical,liang2009multi}:
\begin{equation}
\label{eq:PLV}
P_{LV}(t) = E(t)\left(V(t)-V_0 \right) + S(t) Q_{LV}(t),\;\frac{dV}{dt} = -Q_{LV}(t),
\end{equation}
where $V$ is the volume of the left ventricle and $Q_{LV}$ is the flow rate from the left ventricle into the aorta. $V_0$ is the dead volume of the left ventricle and $S$ 
denotes
the viscoelasticity coefficient of the cardiac wall. For simplicity, we assume that $S$ depends linearly on $P_{LV}$: $S(t) = 5.0 \cdot 10^{-4} \cdot P_{LV}(t)$. 
The time dependent elasticity parameter $E$ is given by \cite{liang2009multi}:
\begin{equation*}
E(t) = E_{max} \cdot e_v(t) + E_{min},
\end{equation*}
\begin{equation*}
\;e_v(t) = 
\begin{cases}
0.5 \left( 1.0-\cos \left( \frac{\pi t}{T_{vcp}} \right) \right),\; 0 \leq t \leq T_{vcp}, \\
0.5 \left( 1.0-\cos \left( \frac{\pi \left( t-T_{vcp} \right)}{T_{vrp}} \right) \right),\; T_{vcp} \leq t \leq T_{vcp}+T_{vrp}, \\
0.0,\;T_{vcp}+T_{vrp} < t \leq T.
\end{cases}
\end{equation*}
$T$ represents the length of the heart cycle. $E_{max}$ and $E_{min}$ are the maximal and minimal elasticity parameters, 
while $T_{vcp}$ and $T_{vrp}$ refer to the durations of the ventricular contraction and relaxation. 
The flow rate $Q_{LV}$ through the aortic valve is governed by the Bernoulli law \cite{sun1997comprehensive} incorporating the viscous resistance and inertia of blood:
\begin{equation}
\label{eq:QLV}
L \frac{dQ_{LV}}{dt} = \Delta P_v(t) - R \cdot Q_{LV}(t) -B \cdot Q_{LV}(t) \cdot \left|Q_{LV}(t) \right|.
\end{equation}
The parameters $R,\;B$ and $L$ quantify the viscous effects, flow separation and inertial effects. Finally, the pressure drop $\Delta P_v(t)$ is computed as follows:
$\Delta P_v(t) = P_{LV}(t)-p_1(0,t),\;t>0$, where $p_1 \left(0,t\right)$ is the pressure at the root of the aorta.

During the systolic phase of the heart cycle, it holds: $\Delta P_v(t)>0$ and we use \eqref{eq:PLV}--\eqref{eq:QLV} to compute
$Q_{LV}(t)$. This value serves as a Dirichlet boundary value at $z_1=0$ for the 1-D model in vessel $V_1$. 
Based on $Q_{LV}(t)$, Equation \eqref{eq:W1} and \eqref{eq:W2} and an approximation of $W_{1,1}\left(0,t \right)$ by \eqref{eq:Extrapolation} the ingoing variable 
$W_{2,1}\left(0,t \right)$ can be determined.

Within the diastolic phase of the heart cycle, $p_1(0,t)$ begins to exceed the pressure in the left ventricle $P_{LV}(t)$. As a result 
the aortic valve is closing and we have no flux or a very little flux between the left ventricle and the aorta and therefore we set $Q_{LV}(t) = 0$.
Since we simulate only the left ventricle without taking into account the filling process by the left atrium, we reactivate
the model at the begin of every heart cycle \cite{formaggia2006numerical}. Thereby, at the end of each heart cycle, the volume of the left ventricle is set to its maximal 
value:
$V \left( k \cdot T \right) = V_{max},\;k \in \mathbb{N}$.

\subsection{Modelling of bifurcations} 

\begin{figure}[h!]
\begin{center}
\includegraphics[scale=0.32]{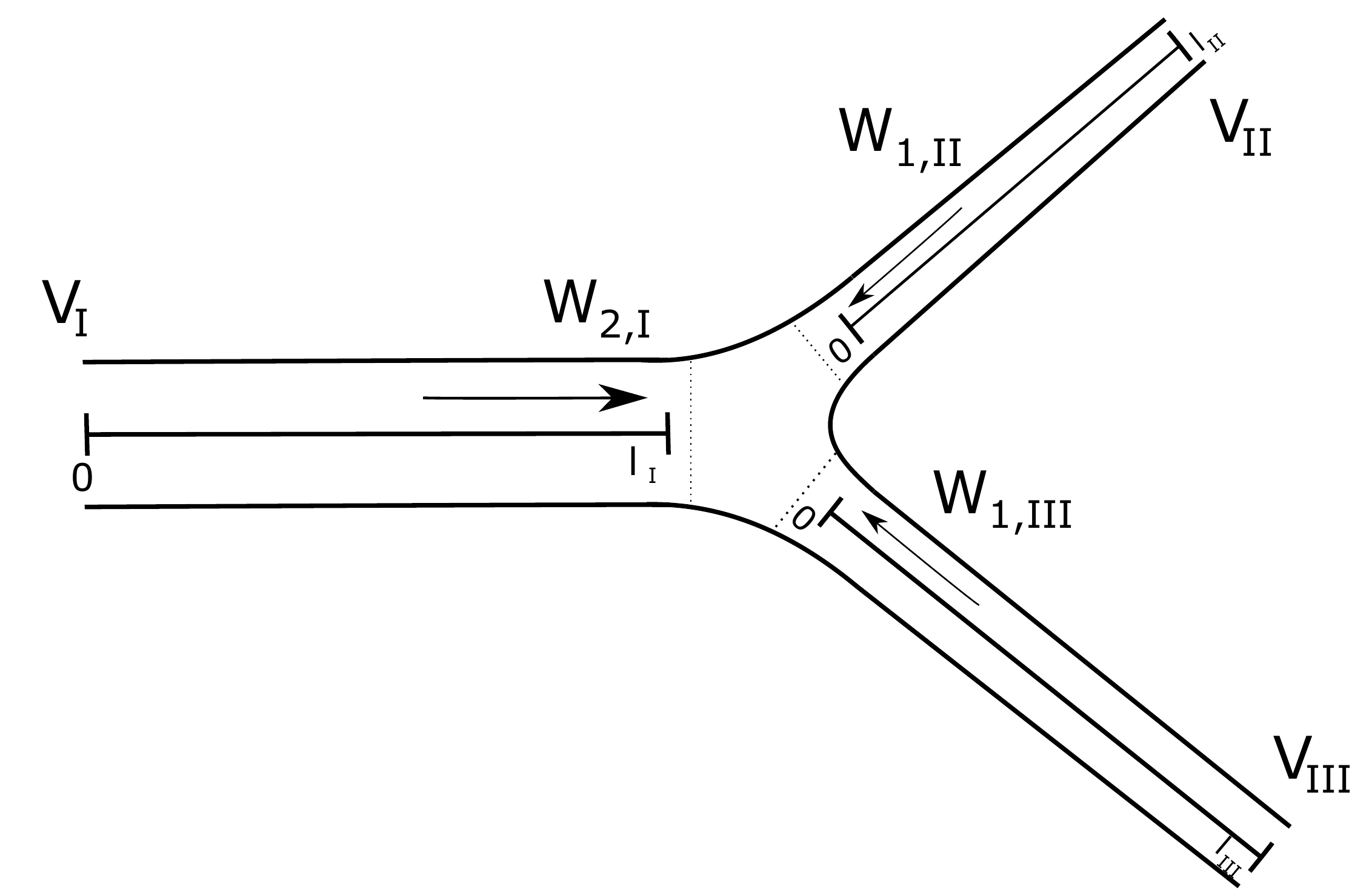}
\end{center}
\caption{\label{fig:Bifurcation} Decomposition of a bifurcation into three independent vessels $V_\text{I}$, $V_\text{II}$ 
and $V_\text{III}$. Orientating the axes of the vessels as in the figure,
the characteristic variables $W_{2,\text{I}}$, $W_{1,\text{II}}$ and $W_{1,\text{III}}$ are leaving the corresponding
vessels. These variables can be determined by values from the inside of the vessels and combined with the coupling conditions \eqref{eq:bif_cond},
we have a system of equations yielding the boundary values for $V_\text{I}$, $V_\text{II}$ and $V_\text{III}$.}
\end{figure}

In order to decrease the flow velocity and to cover the whole body with blood vessels, the arterial system exhibits
several levels of branchings. Therefore, it is very important to simulate blood flow through a bifurcation as exact as possible.
Bifurcations and their mathematical modelling have been the subject of many publications 
\cite{carson2017implicit,formaggia2003one,holden1999riemann,koppl2013reduced,smith2003happens}.
Coupling conditions for systems linked at a bifurcation can be derived by the principles of mass conservation
and continuity of the total pressure. The total pressure for Vessel $V_i$ is defined by:
$$
p_{t,i} = \frac{\rho}{2} \left( \frac{Q_i}{A_i} \right)^2 + p\left(A_i \right).
$$
Indexing the vessels at a bifurcation by $V_i,\;i \in \left\{ \text{I},\text{II},\text{III} \right\} \subset \left\{1,\ldots,55\right\}$, we obtain
the following three coupling conditions:
\begin{equation}
\label{eq:bif_cond}
Q_\text{I} = Q_\text{II} + Q_\text{III} \quad\text{ and }\quad p_{t,\text{I}} = p_{t,\text{II}},\quad p_{t,\text{I}} = p_{t,\text{III}}.
\end{equation}
The remaining equations are obtained by the characteristics entering the bifurcation (see Figure \ref{fig:Bifurcation}). According to Subsection \ref{sec:SingVess}
we have at each bifurcation three characteristics moving from the vessels into the bifurcation. The outgoing characteristic variables can be determined by tracing back
the corresponding characteristic curves (see Subsection \ref{sec:NMC}, Equation \eqref{eq:Extrapolation} and \eqref{eq:WBoundary}). 
Using the characteristic variables and inserting \eqref{eq:W1} and \eqref{eq:W2}
into \eqref{eq:bif_cond}, we obtain a non-linear system of equations for the three unknown ingoing characteristic variables 
$W_{1,\text{I}}$, $W_{2,\text{II}}$ and $W_{2,\text{III}}$. 

\subsection{Modelling of the peripheral circulation}
\label{sec:PerCirc}

At the outlet of a terminal vessel $V_i$,
the reflections of the pulse waves at the subsequent vessels have to be incorporated to simulate a realistic pressure decay. 
For this purpose we assign to each terminal vessel a reflection 
parameter $R_{p,i} = R_{1,i} + R_{2,i}$, where $R_{1,i}$ is the resistance parameter of $V_i$ and $R_{2,i}$ is the equivalent resistance parameter
for all the vessels which are connected to $V_i$ but not contained in the 1-D network. A third parameter $C_i$ quantifies the compliance of the omitted
vessels and, therefore, it is a measure for the ability of these vessels to store a certain blood volume. These parameters form a triple 
$\left(R_{1,i},C_i,R_{2,i} \right)$ that is referred to as a \emph{three-element Windkessel} model in literature 
\cite{alastruey2008lumped} \cite[Chap. 10]{formaggia2010cardiovascular}. In order to describe the dynamics of a Windkessel model the following ODE has
been derived using averaging techniques and an analogy from electrical science \cite{alastruey2007modelling,alastruey2008lumped,marchandise2009numerical}:
\begin{equation}
\label{eq:Outflow}
p_{i,t} + R_{2,i} C_i \frac{d p_{i,t}}{dt} = p_v + \left( R_{1,i} + R_{2,i} \right) Q_{i,t} + R_{1,i} R_{2,i} C_i \frac{d Q_{i,t}}{dt}.
\end{equation}
$p_{i,t} = p\left(A_i\right)$ and $Q_{i,t}$ denote the pressure and flow rate at the outlet $z=l_i$ of a terminal vessel $V_i$, respectively.
$p_v$ is an averaged pressure in the venous system. 
Combining \eqref{eq:Outflow} with \eqref{eq:PressureAreaRelation}, \eqref{eq:W1} and \eqref{eq:W2} yields an equation depending on the characteristic variables. 
By means of this equation 
and the given outgoing characteristic variable $W_{2,i}$, the missing ingoing characteristic variable $W_{1,i}$ can be computed, by solving for each time point of
interest a non-linear equation. 
Having $W_{1,i}$ and $W_{2,i}$ at $z=l_i$ and for a time point $t>0$ at hand, 
the boundary values $A_i\left(l_i,t\right)$ and $Q_i\left(l_i,t\right)$ can be computed for
each $t>0$ using \eqref{eq:W1} and \eqref{eq:W2}. Further information on the derivation of lumped parameter models for the peripheral circulation 
can be found in \cite{olufsen2004deriving}.

\subsection{Modelling the influence of a stenosis on blood flow}
\label{sec:Stenosis}

\begin{figure}[h!]
\begin{center}
\includegraphics[scale=0.25]{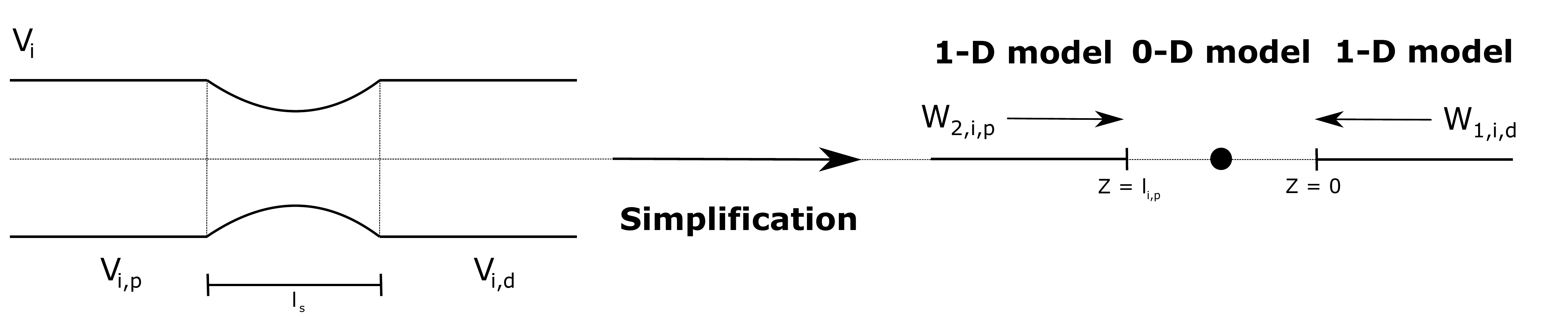}
\caption{\label{fig:Stenosis_0D_Model} Decomposition of a blood vessel $V_i$ into a proximal part $V_{i,p}$ and a distal part $V_{i,d}$. The stenosis having
a degree of stenosis $R_s \in \left[0,1 \right]$ is lumped to a node and modelled by a $0$-D model. In order to provide boundary conditions at $z=l_{i,p}$ and
$z=0$ for the $1$-D models accounting for the influence of the stenosis, the outgoing characteristics $W_{2,i,p}$ and $W_{1,i,d}$ are combined with the coupling
conditions of the $0$-D model.}
\end{center}
\end{figure}

A blood vessel $V_i$ containing a stenosis is split into three parts: A proximal part $V_{i,p}$, the stenosis itself and a distal part $V_{i,d}$.
In a next step, we assign to $V_{i,p}$ and $V_{i,d}$ the 1-D blood flow model from Subsection \ref{sec:SingVess}, while the part of $V_i$ that is
covered by the stenosis is lumped to a node (see Figure \ref{fig:Stenosis_0D_Model}). 
The degree of stenosis is represented by a parameter $R_s \in \left[ 0,1 \right]$, where $R_s = 0$ corresponds to 
the healthy state and $R_s = 1$ stands for the case of a completely occluded blood vessel. For convenience, it is assumed that both vessel parts
have the same section area $A_{0,i}$ as well as the same elasticity parameters. The lengths of $V_{i,p}$ and $V_{i,d}$ are denoted by $l_{i,p}$ and
$l_{i,d}$. At the boundaries $z_{i,p} = l_{i,p}$ and $z_{i,d} = 0$ that are adjacent to the stenosis, two characteristic variables $W_{2,i,p}$ and $W_{1,i,d}$
are moving towards the stenosis (see Figure \ref{fig:Stenosis_0D_Model} and Subsection \ref{sec:SingVess}). Using these characteristic variables
and an appropriate $0$-D model, we have enough equations to compute the boundary conditions for the two 1-D models. Modelling the stenosis, we consider 
an ODE containing several parameters of physical relevance to couple both parts of the affected vessels. 
According to \cite{seeley1976effect,stergiopulos1992computer,young1973flow} the flow rate $Q_s$ through a stenosis and 
the pressure drop $\Delta p_s\left(t\right) = p_{i,d}\left(0,t\right)-p_{i,p}\left(l_{i,p},t\right)$ across a stenosis are related to each 
other by the following ODE:
\begin{equation}
\label{eq:ODESten}
\frac{K_u \cdot \rho \cdot l_s}{A_{0,s}} \frac{dQ_s}{dt} = \Delta p_s - \frac{K_v \cdot \eta \cdot l_s}{A_{0,s} \cdot D_s}Q_s 
- \frac{K_t \cdot \rho}{2A_{0,s}^2} \left( \frac{A_{0,s}}{A_s}-1 \right)^2 Q_s \left|Q_s \right|.
\end{equation}
$A_{0,s}$ and $A_s = \left(1-R_s \right) \cdot A_{0,s}$ refer to the section areas of the normal and stenotic segments, while $D_0$ and $D_s$ denote the 
corresponding diameters \cite{liang2009multi}. $l_s$ is the length of the stenosis. The remaining parameters are empirical coefficients, which are given by 
\cite{liang2009multi}:
\begin{equation*}
K_v = 32.0 \cdot \left( 0.83 \cdot l_s + 1.64 \cdot D_s \right) \cdot\left( \frac{A_{0,s}}{A_s} \right)^2 \frac{1}{D_s},\;K_t = 1.52 \;\text{ and }\; K_u = 1.2. 
\end{equation*}
Solving \eqref{eq:ODESten}, we use the solution value $Q_s$ for each time point as a boundary condition at $z_{i,p} = l_{i,p}$ and $z_{i,d} = 0$. 
Together with the extrapolated characteristic variables
$W_{2,i,p}$ and $W_{1,i,d}$ as well as \eqref{eq:W1} and \eqref{eq:W2}, we have a system of equations for the boundary conditions adjacent to the stenosis. 

In the case of a full occlusion, i.e. $R_s = 1$, we multiply \eqref{eq:ODESten} by $A_s^2$. Considering the limit $A_s \rightarrow 0$, we obtain an algebraic 
equation that gives $Q_s = 0$. 
This is equivalent to a full reflection of the ingoing characteristics $W_{2,i,p}\left( l_{i,p},t \right)$ and $W_{1,i,d}\left( 0,t \right),\;\;t>0$.

\section{Kernel-based surrogate models}
\label{sec:KerMeth}
In this section we introduce kernel methods for surrogate modeling. First, we present the general ideas of kernel methods applied to the 
approximation of an arbitrary continuous function $f: \Omega\subset\mathbb R^d\to \mathbb R^q$, where $\Omega$ is a given input parameter domain and $d, q$ are the input 
and 
output dimensions, which can be potentially large (e.g., the present setting will lead to $q=400$). Then we concentrate on the present field of application and 
discuss how this general method can be used to 
produce cheap surrogates of the full model described in the previous sections.

Our goal is to construct a surrogate function ${\hat f} : \Omega \to \mathbb R^q$ such that ${\hat f} \approx f$ on $\Omega$, while the evaluation of ${\hat f}$ 
for 
any input value is considerably cheaper than evaluating $f$ for the same input. 
This approximation is produced in a data-dependent fashion, i.e., a finite set of snapshots, obtained with the full simulation, is used to train the model to 
provide a good prediction of the exact result for any possible input in $\Omega$. The computationally demanding construction of the snapshots is performed only 
once and in 
an offline phase, while the 
online computation of the prediction for a new input parameter uses the cheap surrogate model.

\subsection{Basic concepts of kernel methods}
We introduce here only the basic tools needed for our analysis. For an extensive treatment of kernel-based approximation methods we refer to the monographs \cite{ 
Fasshauer2007, Fasshauer2015, Wendland2005}, while a  detailed discussion of kernel-based sparse surrogate models can be found in \cite{HS2017a}. Nevertheless, 
we recall that this technique has several advantages over other approximation methods, namely it allows for large input and output dimensions, it works with 
scattered data, it allows fast and sparse solutions through greedy methods and has a notable flexibility related to the choice of the particular kernel.

As recalled before, we aim at the reconstruction of a function $f:\Omega\to \mathbb R^q$, $\Omega\subset \mathbb R^d$. We assume to 
have a dataset given by $N\in\mathbb N$ pairwise distinct inputs, i.e., a set $X_N:=\{x_1, \dots, x_N\}\subset \Omega$ of $N$ points in $\Omega$ (the data 
points) and corresponding function 
evaluations $F_N:=\{f(x_i), x_i\in X_N\}\subset\mathbb R^q$ (the data values). 

The construction of ${\hat f}$ makes use of a positive definite kernel $K$ on $\Omega$. We recall that a function $K:\Omega\times \Omega\to \mathbb R$ is a 
strictly 
positive kernel on $\Omega$ if it is symmetric and, for any $N\in\mathbb N$ and any set of pairwise distinct points $X_N:=\{x_1, \dots, 
x_N\}\subset \Omega$, the $N\times N$ kernel matrix $A_{K, X_N} := (K(x_i, x_j))_{i,j=1}^N$ is positive definite. Many strictly positive definite kernels are 
known in explicit 
form, and 
notable examples are e.g. the Gaussian $K(x, y):= \exp(-\varepsilon^2\|x-y\|_2^2)$ (where $\varepsilon$ is a tunable parameter) and the Wendland kernels 
\cite{Wendland1995a}, which are radial and compactly supported kernels of 
piecewise polynomial type and of finite smoothness. Given a kernel $K$, the surrogate model ${\hat f}: \Omega \to \mathbb R^q$ is constructed via the ansatz 
\begin{equation}\label{eq:kernel_ansatz}
 {\hat f}(x) := \sum_{j=1}^N \alpha_j K(x, x_j), \; x\in\Omega,
\end{equation}
with unknown coefficient vectors $\alpha_j\in\mathbb R^q$. The coefficients are obtained by the vectorial interpolation conditions  
\begin{equation}\label{eq:kernel_int_conditions}
 {\hat f}(x_i) = f(x_i), \; 1\leq i\leq N,
\end{equation}
i.e., the surrogate model ${\hat f}$ is required to predict the same value of the full model $f$ when computed on each of the data points contained in the 
dataset.
Putting together the ansatz \eqref{eq:kernel_ansatz} and the interpolation conditions \eqref{eq:kernel_int_conditions}, one obtains the set of equations
\begin{equation*}
 {\hat f}(x_i) = \sum_{j=1}^N \alpha_j K(x_i, x_j)= f(x_i), \; 1\leq i\leq N,
\end{equation*}
which can be formulated as a linear system $A_{K, X_N}\alpha = b$ with 
\begin{equation}\label{eq:matrices}
A_{K, X_N} := [K(x_i, x_j)]_{i,j=1}^N \in \mathbb R^{N\times N},\; \alpha := \left[\begin{array}{c}
                                                                        \vdots\\
                                                                        \alpha_j^T\\
                                                                        \vdots\\
                                                                        \end{array}\right]\in\mathbb R^{N\times q}, 
                                                                        \;b := \left[\begin{array}{c}
                                                                        \vdots\\
                                                                        f(x_i)^T\\
                                                                        \vdots\\
                                                                        \end{array}\right]\in\mathbb R^{N\times q}.
\end{equation}

Since the kernel is chosen to be strictly 
positive definite, the matrix $A_{K, X_N}$ is positive definite for any $X_N$, thus the above linear system possesses a unique solution $\alpha \in \mathbb 
R^{N\times q}$. In 
other terms, the model \eqref{eq:kernel_ansatz} satisfying interpolation conditions \eqref{eq:kernel_int_conditions} is uniquely defined for arbitrary pairwise 
distinct data 
points $X_N$ and data values $F_N$.

This interpolation scheme can be generalised by introducing a regularisation term, which reduces possible oscillations in the surrogate at the price of a non 
exact 
interpolation of the data. We remark that in principle this does not reduce the accuracy, since a parameter $\lambda\geq0$ can be used to tune the influence of 
the regularization 
term, and a 
zero value can be used when no regularization is needed. 

To explain this in details, we first recall that, associated with a strictly positive definite kernel there is a uniquely defined Hilbert space $\mathcal 
H_K(\Omega)$ of functions 
from $\Omega$ to $\mathbb R^q$. For the sake of simplicity, we discuss the case $q=1$, i.e., scalar valued functions, while the generalisation to vectorial 
functions will be 
sketched at the end of this section. The space $\ns$ contains in particular all the functions of the form \eqref{eq:kernel_ansatz}, and their squared norm can 
be computed as $\|\hat f\|^2_{\ns} = \alpha^T A_{K, X_N} \alpha$. This means that a surrogate with small $\ns$-norm is defined by coefficients $\alpha$ 
with small magnitude.

With these tools, and again for $q=1$ and a regularisation parameter $\lambda\geq0$, a different surrogate can be defined as the solution of the optimisation 
problem
\begin{equation}\label{eq:representer}
\min_{{\hat f}\in\ns} \sum_{i=1}^N \left(f(x_i) - {\hat f}(x_i)\right)^2 + \lambda \|{\hat f}\|_{\ns}^2 = \min_{\alpha\in\mathbb{R}^N} \|A_{K, X_N} \alpha - 
b\|_2^2 + \lambda \alpha^T A_{K, X_N} \alpha,
\end{equation}
which is a regularised version of the interpolation conditions \eqref{eq:kernel_int_conditions}, where exact interpolation is 
replaced by square error minimization, and the 
surrogate is requested to have a small norm. When a strictly positive definite kernel is used, the Representer Theorem \cite{Schoelkopf2001v} guarantees that 
the problem 
\eqref{eq:representer} has a unique solution, that this solution is of the form \eqref{eq:kernel_ansatz}, and that the coefficients $\alpha$ are defined as the 
solutions of the 
linear system 
\begin{equation*}
 \left(A_{K, X_N} + \lambda I\right) \alpha = b,
\end{equation*}
where $I$ is the $N\times N$ identity matrix, and where now $\alpha$ and $b$ are column vectors. It is now clear that pure interpolation can be obtained by 
letting $\lambda = 
0$ although a positive $\lambda$ 
improves the conditioning of the linear system, reducing possible oscillations in the solution.

These Hilbert spaces, the corresponding error analysis, and the formulation of the regularised interpolant can be extended to the case of vector valued 
functions 
$f:\Omega\to\mathbb R^q$ just by applying the same theory to each of the $q$ components. 
The fundamental point here, to have an effective method to be used in surrogate modeling, is to avoid having $q$ different surrogates, one for each component. 
This 
would result in $q$ independent sets of centers, hence many kernel evaluations to compute a point value ${\hat f}(x)$. To reduce the overall number of centers, 
one can make the 
further assumption that a common set of centers is used for all components. 
From the point of view of the actual computation of the interpolant, this is precisely equivalent to the solution of the linear system \eqref{eq:matrices}, 
where in the regularised 
case also the term $\lambda I$ is included. We remark that more sophisticated approaches are possible to treat vector valued functions, but the approach 
presented in this work
yields already satisfactory results.

\subsection{Sparse approximation}
So far, we have shown that kernel interpolation is well defined for arbitrary data, and that the corresponding interpolant has certain approximation properties.
Although the method can deal with arbitrary pairwise distinct inputs $X_N$, the resulting surrogate model ${\hat f}$ is required to be fast to evaluate. From 
formula \eqref{eq:kernel_ansatz}, it is clear that the computational cost of the evaluation of ${\hat f}(x)$ on a new input parameter $x\in \Omega$ is 
essentially related to the 
number $N$ 
of elements in the sum. It is thus desirable to find a sparse expansion of the form \eqref{eq:kernel_ansatz}, 
i.e., one where most of the coefficient vectors $\alpha_j$ are zero. This sparsity structure can be obtained by selecting a small subset $X_n\subset X_N$ of the 
data points, and 
computing the corresponding surrogate. 

An optimal selection of these points is a combinatorial problem, which is too expensive with respect to the computational effort. 
Instead, we employ greedy methods (see \cite{Temlyakov2008} for a general treatment, and \cite{DeMarchi2005, SchWen2000} for the case of kernel approximation). 
Such methods select 
a sequence of data points starting with the empty set $ X_0:=\emptyset$, and, at iteration $n\geq 1$, they update the set as $X_n:=X_{n-1}\cup\{x_n\}$ by adding 
a suitable 
selected point $x_n\in X_N\setminus X_{n-1}$. The selection of $x_n$ is done here with the $f$-Vectorial Kernel Orthogonal Greedy Algorithm ($f$-VKOGA, 
\cite{Wirtz2013}), which 
works as follows.
At each iteration, a partial surrogate can be constructed as 
$$
{\hat f}_n(x) = \sum_{j=1}^n \alpha_j K(x, x_j),\;\; (A_{K, X_n} + \lambda I)\alpha = b_{X_n}, 
$$
where $A_{K, X_n}$, $b_{X_n}$ are the matrix and vector of \eqref{eq:matrices} restricted to the points $X_n$. To evaluate the quality of the partial 
surrogate, one can check the 
residual vector
$$
r_n(f)(x) := f(x) - {\hat f}_n(x)
$$
for all $x\in X_N\setminus X_n$. The $f$-VKOGA takes precisely $x_n:= \max_{x\in X_N \setminus X_{n-1}} \|r_n(x)\|_2$, i.e., it includes in the model the data 
point where the 
error is currently largest. By checking the size of the residual, one can stop the iteration with ${\hat f}_n\approx {\hat f}$, while potentially $n\ll N$, 
i.e., 
the new surrogate 
model is much cheaper to evaluate but it retains the same accuracy of ${\hat f}$. More precisely, it has been proven that, under smoothness assumptions on the 
target function $f$, the 
VKOGA algorithm, with the $f$- or similar selection strategies, can attain algebraic or even spectral convergence rates \cite{SH16b, Wirtz2013}.

Finally, we remark that the partial surrogates can be efficiently updated when adding a new point, i.e., ${\hat f}_n$ can be obtained from $\hat f_{n-1}$ by 
computing only a new 
coefficient in the expansion, while the already computed ones are not modified. We point to the paper \cite{HS2017a} for a more in depth explanation of this 
efficient computational process.

Observe that this greedy method results in the selection of a small subspace $V_n:=\mathrm{span}\{K(\cdot, x_i), x_i\in X_n\}$, and $\hat f_n$ is 
computed as the projection, thus best approximation, of $f$ into $V_n$. The selection of the points $X_n$ via the $f$-greedy selection strategy makes use of the values 
of 
$f$ on all the points $X_N$. In this sense, the procedure is similar to a least square approximation, where a small set of points is used to generate an accordingly 
small approximation space. Nevertheless, it is not clear in the least square setting how these few points should be selected, whereas the present approach allows 
an incremental selection of points and an efficient update of the approximant, which can be stopped when a tolerance criterion is reached. Moreover, by solving equation 
\eqref{eq:representer} we are indeed constructing an approximant that minimizes a least squares accuracy term combined with a regularization term.

\subsection{Simulating blood flow in the vicinity of a peripheral stenosis by means of kernel methods}
Coming back to the blood flow simulation, we define in details the target functions $f$ which will be approximated by the kernel method. These functions will 
represent the maps 
from an input stenosis degree $R_s\in[0,1]$ to the resulting pressure or flow-rate curve for different vessels, as computed by the full simulation of Section 
\ref{sec:NMC}. The definition of $f$ is described in the following.

Since the numerical simulation is expected to have a transient phase before reaching an almost-periodic behaviour, the system is first simulated with the 
method 
of Section \ref{sec:NMC} in the time interval $[T_0:=0s, T_1:=20s]$ for the healthy state $R_s = 0$. The state reached at time $T_1$ is then used as initial 
value for the 
subsequent simulations. At the time $T_1$ the stenosis is activated with a degree $R_s\in[0,1]$ and the system is simulated until $T_2=30s$ for 
various values of $R_s$. 

From this set of simulations for different values of $R_s$, we keep the pressure and flow-rate curves of the last heart beat, i.e., in the time interval $[29s, 
30s]$. 
This means that for each point in the spatial grid we have the time evolution of the pressure and flow rate, which are represented as a $q$-dimensional 
vector for each 
space point, where $q$ depends on the actual time discretization step.

In order to study the effect of the stenosis, we concentrate on the vessels number $N_v \in \{52, 54, 55, 56\}$, which are the ones surrounding the stenosis 
(see Figure 
\ref{fig:ArterialNetwork_Stenosis}), and for each of those we select a reference space point.
Putting all together, for each of the four vessels we have one reference point located in the middle
of these vessels. For each point we have the $q$-dimensional time discretization of the pressure and flow rate 
curve in the time interval $[29s, 30s]$. The maps from an input stenosis degree to these vectors define functions $f^P_{N_v}: [0,1]\to\R^q$ (for the pressure) 
and $f^{F}_{N_v}: 
[0,1]\to\R^q$ (for the flow rate). 

\begin{figure}[h!]
\begin{center}
\includegraphics[scale=0.4]{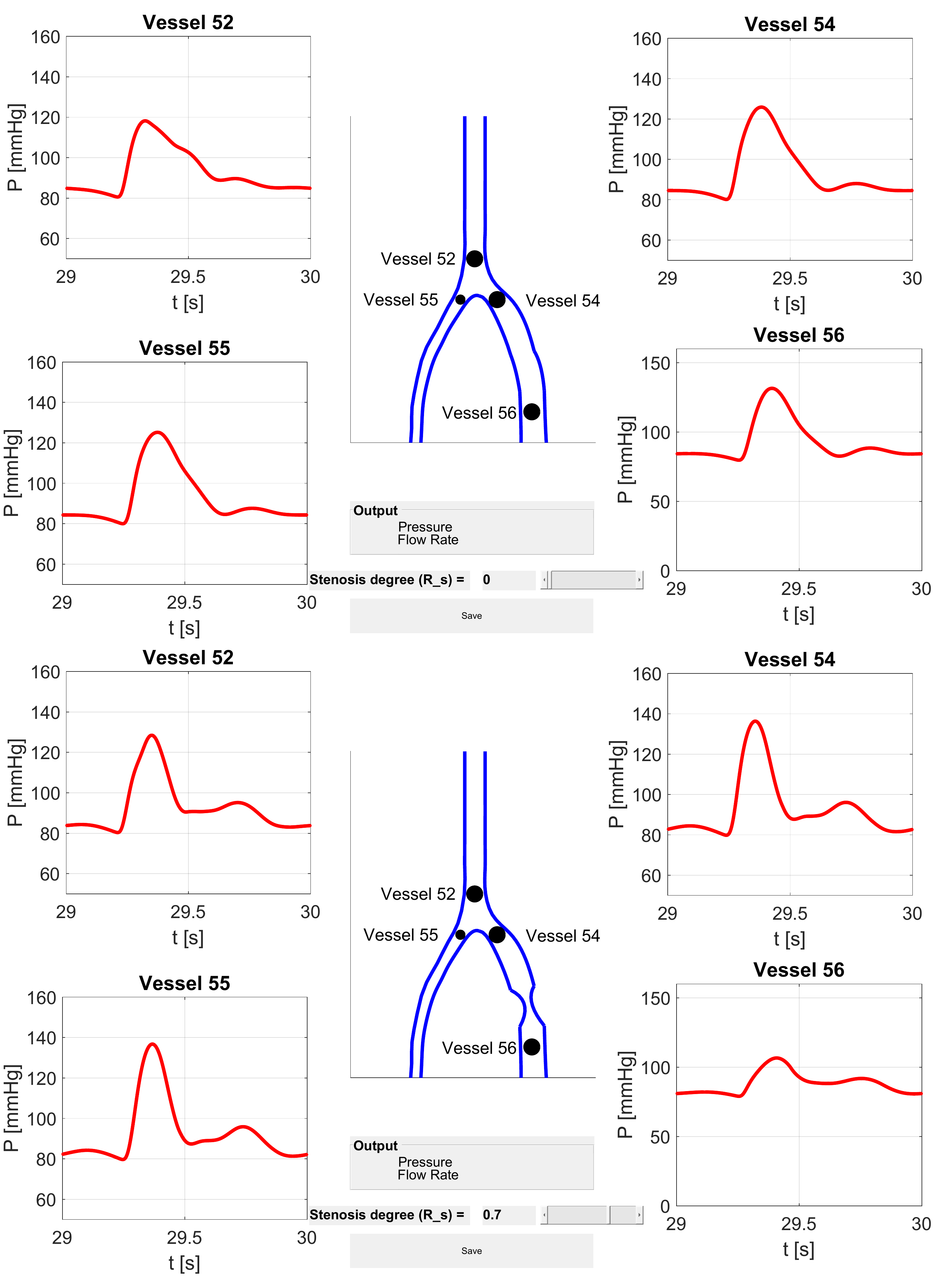}
\caption{\label{fig:example_pressure_curves} Pressure curves around the stenosis for the healthy state (top) and a degree of stenosis $R_s \approx 0.7$ 
(bottom).
The curves are reported at the black dots for $t\in \left[29\unit{s},30\unit{s} \right]$.}
\end{center}
\end{figure}

\begin{figure}[h!]
\begin{center}
\includegraphics[scale=0.4]{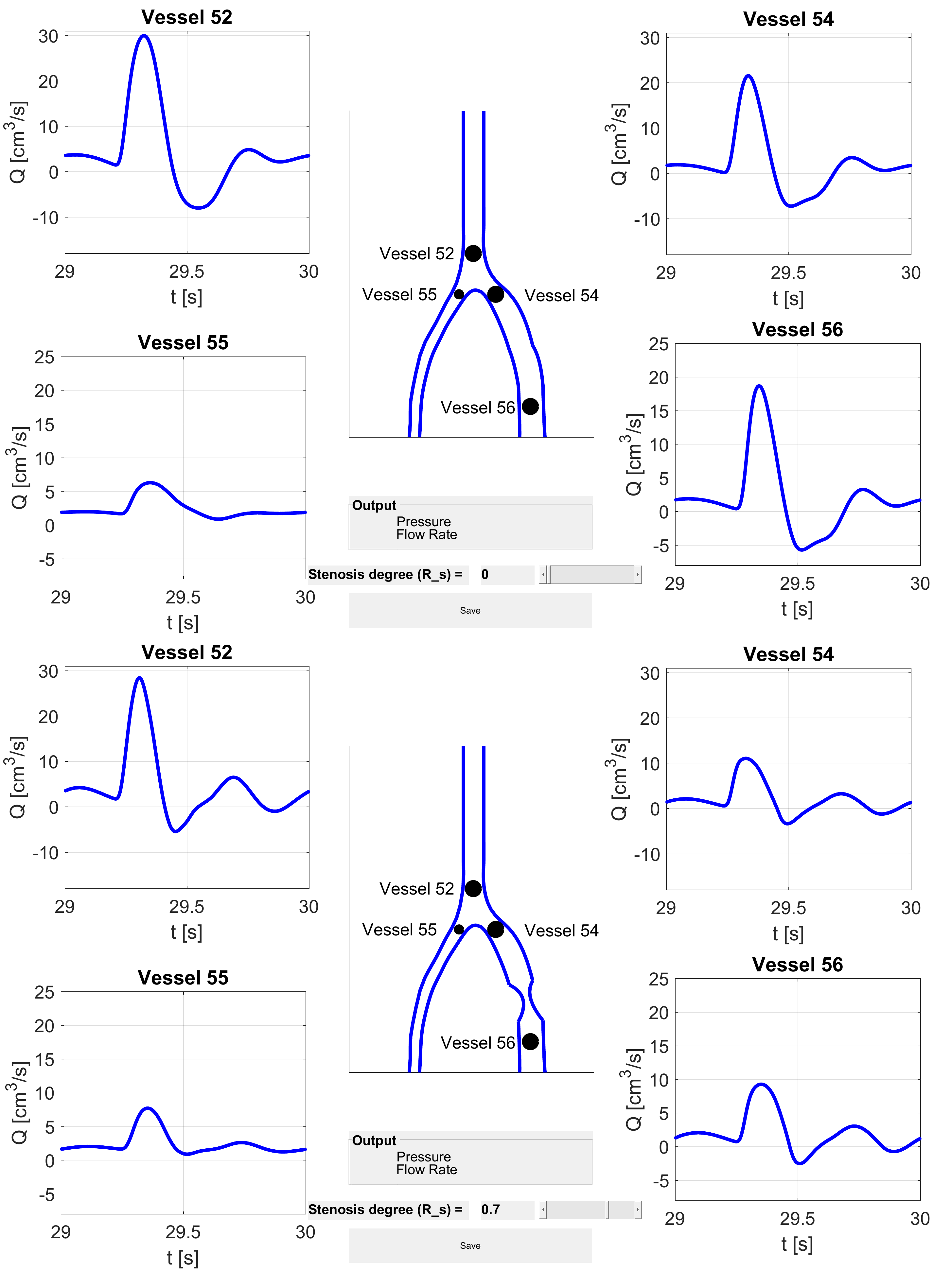}
\caption{\label{fig:example_flow_rate_curves} Flow rate curves around the stenosis for the healthy state (top) and a degree of stenosis $R_s \approx 0.7$ 
(bottom).
The curves are reported at the black dots for $t\in \left[29\unit{s},30\unit{s} \right]$.}
\end{center}
\end{figure}

Figure \ref{fig:example_pressure_curves} and \ref{fig:example_flow_rate_curves}
show examples of pressure and flow rate curves, for both healthy state and for $R_s \approx 0.7$. It can be observed that in the case $R_s \approx 0.7$ 
the flow rate in Vessel $54$ and $56$ is remarkably reduced, while the flow rate in Vessel $55$ is slightly enlarged.
The pressures in Vessel $52$, $54$ and $55$ are increased, which may lead to the formation of an aneursysm, if the vessel walls are weakend in this region. 
Concerning the healthy state $R_s = 0$, one has to note that the pressure values are within a physiological reasonable range, i.e., $79 \unit{mmHg}$
for the diastolic pressure and $130 \unit{mmHg}$ for the systolic pressure. Similar pressure curves have also been 
published in other works \cite{alastruey2008lumped,liang2009multi}.

With respect to the general setting introduced in the previous section, we have here $d:=1$, $q=400$, $\Omega := [0,1]$. We can then train a kernel model for 
each of the eight functions 
(corresponding to pressure and flow-rate for each of the four vessels) using 
snapshot-based datasets. In particular, the data points $X_N\subset[0,1]$ are a set of $N$ pairwise distinct stenosis degrees, and $F_N\subset\mathbb R^{q}$ is 
the 
set of snapshots obtained by the full model run on the input parameters $X_N$. 

Those models can then be used to predict the output of the simulation for an input stenosis not present in the dataset. For example, for a given value $R_s\in 
[0,1]$, the 
evaluation $\hat f_{N_v}^{P}(R_s)$ is a $q$-dimensional vector which approximates the pressure curve in the time interval $[29s, 30s]$ in vessel $N_v$ with 
stenosis degree $R_s$. We remark that this setting can be easily modified to approximate different aspects of the full simulation, although the present ones 
yield interesting insights into the behaviour of the system.

Before we present the numerical tests, two remarks on the data are in order. First, the current time step produces $400$ samples per second, which means that 
we have $q=400$, i.e., we are approximating functions $\left[0,1 \right] \to \R^{400}$.
Second, the data obtained with $R_s=0$ are removed from the datasets and replaced with the one with $R_s = 10^{-6}$, since the ODE model \eqref{eq:ODESten} is meaningful 
only for a strictly positive value of the stenosis degree, while the quadratic term vanishes for $R_s=0$, thus leading to a different model.
This restriction is not relevant from an application point of view, since the a value $R_s = 10^{-6}$ can be considered to effectively represent the healthy state.

In the training of each model, 
the parameters $\varepsilon$ and $\lambda$ are chosen within a range of possible values by $k$-fold cross validation. This means that the training 
data are randomly permuted and divided into $k$ disjoint subsets of approximately the same size and, for each pair $(\varepsilon, \lambda)$ of possible parameters, a 
model is 
trained on the union of $k-1$ subsets and tested on the remaining one. This operation is repeated $k$ times changing in all possible ways the $k-1$ sets used for 
training. The 
average of the error obtained by these $k$ tests is assigned as the error score of the parameter pair, and the best pair $(\varepsilon, \lambda)$ is chosen as the one 
yielding the 
smallest error score. The actual model is then trained on the whole training set using these parameters. In more details, we use here $10$-fold cross 
validation and test $20$ 
logarithmic equally spaced values $\varepsilon \in [10^{-2}, 50]$ and $15$ logarithmic equally spaced 
values $\lambda \in [10^{-16}, 10^{-2}]$. The error measure to sort the parameters is the maximum absolute error.

\section{Numerical tests}\label{sec:NumTests}
This section is concerned with the training of the surrogate models and analyses them under different perspectives.
We train different surrogates using different training sets of increasing size to analyse the number of full model runs needed to have an accurate surrogate. 
They are obtained each with $N$ equally spaced stenosis degrees in $\left[10^{-6},1 \right]$, with $N = \{5, 10, 20, 40, 80, 160, 320, 640\}$.

A further dataset of $N = 1000$ equally spaced stenosis degrees is used as a test set, i.e., the surrogates computed with the various training datasets are 
evaluated on this set of input 
stenosis degrees, and the results are compared with the full model computations. We remark that the values of this test set are not contained in the training 
sets (except for 
$R_s=10^{-6}$ and $R_s=1$), so the results are reliable assessments of the models' accuracy. For every value $(R_s)_i$ in the test set we consider the absolute 
and relative errors 
$$
e_A^{(i)}:=\|f((R_s)_i) - {\hat f}((R_s)_i)\|_2,\;\; e_R^{(i)}:=\frac{\|f((R_s)_i) - {\hat f}((R_s)_i)\|_2}{\|f((R_s)_i)\|_2},
$$
and, to measure the overall error over the test set, we compute both the maximum absolute and relative error, i.e., 
$$
E_{A}:= \max_{1\leq i\leq 1000} e_A^{(i)}, \;\;E_{R}:= \max_{1\leq i\leq 1000} e_R^{(i)}. %E_{RMS}:= \frac{1}{\sqrt{1000}}\sqrt{\sum_{i= 1}^{1000} e_i^2}.
$$
We use the Gaussian kernel, and the $f$-VKOGA is stopped using a tolerance $5 \cdot 10^{-8}$ on the regularized \textit{Power Function}, which controls the 
model stability \cite{HS2017a}. 

\subsection{Simulation parameters}

Before we study the performance of the numerical model, we summarise the simulation parameters in this subsection. For the 1-D arterial network, the data
from \cite{sherwin2003computational}[Tab. 1] have been used. In this table the different lengths $l_i$, section areas $A_{0,i}$ and elasticity parameters $\beta_i$ can 
be found. By means of $\beta_i$ and $A_{0,i}$, the elasticity parameters
$G_{0,i}$ in \eqref{eq:PressureAreaRelation} can be calculated as follows: $G_{0,i} = \beta_i \cdot \sqrt{A_{0,i}}$.
Please note that Vessel $54$ from the original data set is split into a new Vessel $54$ of length $l_{54} = 10.0\;\unit{cm}$,
the stenosis of length $l_s = 1.0\;\unit{cm}$ and an additional Vessel $56$ of length $l_{56} = 21.2\;\unit{cm}$ (see Figure \ref{fig:ArterialNetwork_Stenosis}
and Figure \ref{fig:Stenosis_0D_Model}). The resistances $R_{p,i} = R_{1,i} + R_{2,i}$ and capacities $C_i$ occuring in
\eqref{eq:Outflow} are listed in \cite{stergiopulos1992computer}[Tab. 2],
where the resistance $R_{1,i}$ is determined by the characteristic impedance $Z_i = \rho \cdot c\left(A_{0,i} \right)/A_{0,i}$ \cite{alastruey2008lumped}. 
For $t=0$, we set $A_i(z,0)= A_{0,i},\;Q_i(z,0)= 0$ and $p_i(z,0)= 0$ in the corresponding vessels. The parameters for the left ventricle model
are listed in Table \ref{tab:LeftVentricle}. In order to solve the ODEs occuring in Section \ref{sec:DimRedModel}, we use an explicit discretisation of first order.
\begin{table}[h!]
\begin{center}
\caption{\label{tab:LeftVentricle} 
List of the different parameters for the model of the left ventricle \cite{formaggia2006numerical}[Section 4],\cite{liang2009multi}[Tab. 2].}
\begin{tabular}{|c|c|c|c|}
 \hline
 \hline
 Physical Parameter & sign & value & unit\\
 \hline
 dead volume left ventricle & $V_0$ & $10.00$   &  $\unit{cm^3}$  \\
 \hline
 maximal volume    & $V_{max}$  & $130.00$  &  $\unit{cm^3}$  \\
 \hline
 duration of heart cycle & $T$ & $1.00$  &  $\unit{s}$  \\
 \hline
 duration of ventricular contraction & $T_{vcp}$ & $0.30$  &  $\unit{s}$  \\
 \hline
 duration of ventricular relaxation & $T_{vrp}$ & $0.15$  &  $\unit{s}$  \\
 \hline
 maximal elastance & $E_{max}$ & $2.75$  &  $\unit{mmHg/cm^3}$  \\
 \hline
 minimal elastance & $E_{min}$ & $0.08$  &  $\unit{mmHg/cm^3}$  \\
 \hline
 viscous resistance & $R$ & $3.0 \cdot 10^{-3}$  &  $\unit{mmHg\cdot s/cm^3}$  \\
 \hline
 separation coefficient & $B$ & $2.5 \cdot 10^{-5}$  &  $\unit{mmHg\cdot s^2/cm^6}$  \\
 \hline
 inductance coefficient & $L$ & $5.0 \cdot 10^{-4}$  &  $\unit{mmHg\cdot s^2/cm^3}$  \\
\hline
\hline
\end{tabular}
\end{center}
\end{table}

\subsection{Accuracy of the surrogate models}

\begin{figure}[ht]
\begin{center}
\begin{tabular}{cc}
\includegraphics[scale=0.38]{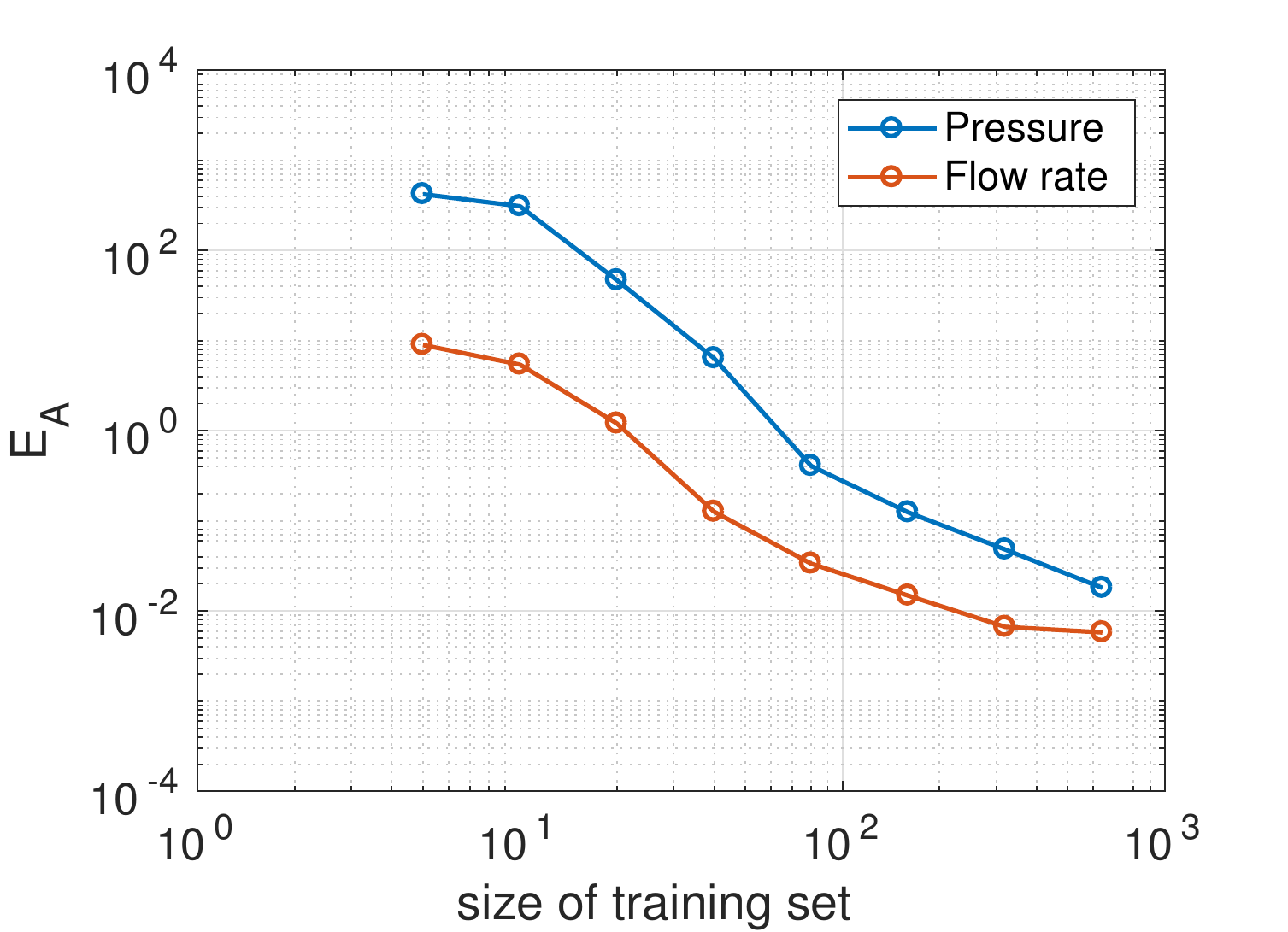}&
\includegraphics[scale=0.38]{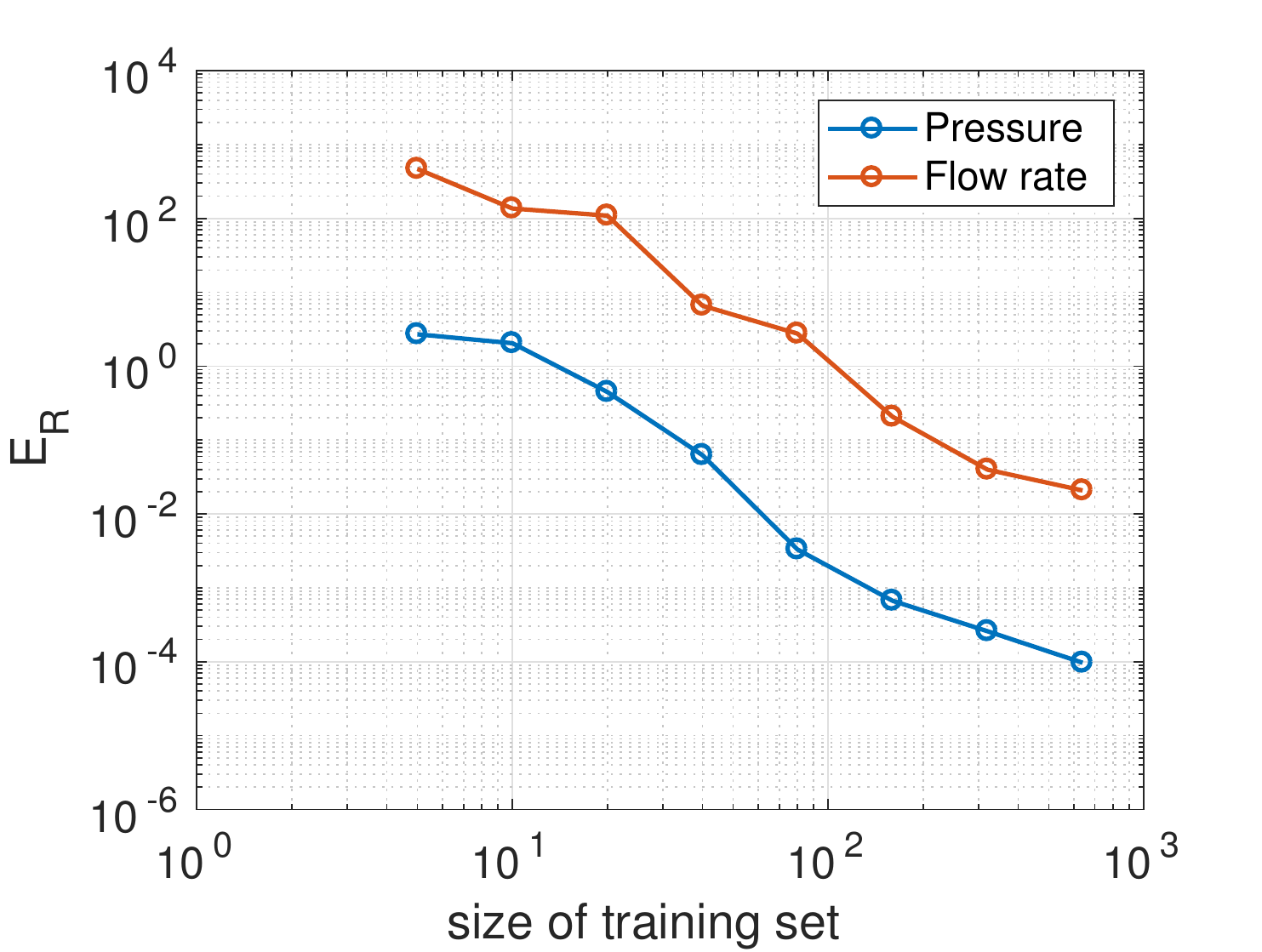}
\end{tabular}
\caption{Absolute (left) and relative (right) errors obtained by kernel surrogates with an increasing training set, for both the 
pressure (blue curves) and flow rate (red 
curves).}\label{fig:errors56}
\end{center}
\end{figure}

\begin{figure}[ht]
\begin{center}
\begin{tabular}{cc}
\includegraphics[scale=0.4]{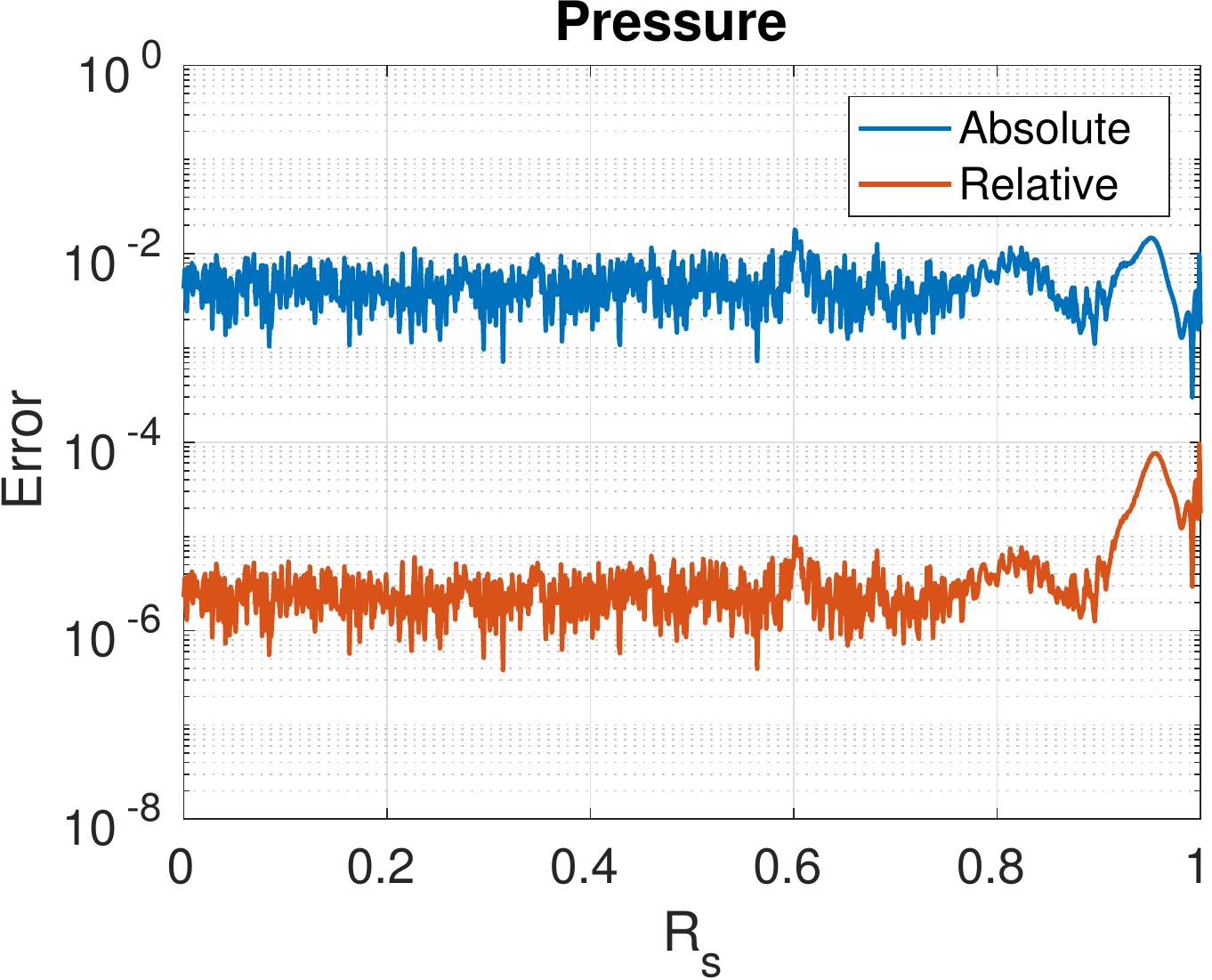}&
\includegraphics[scale=0.4]{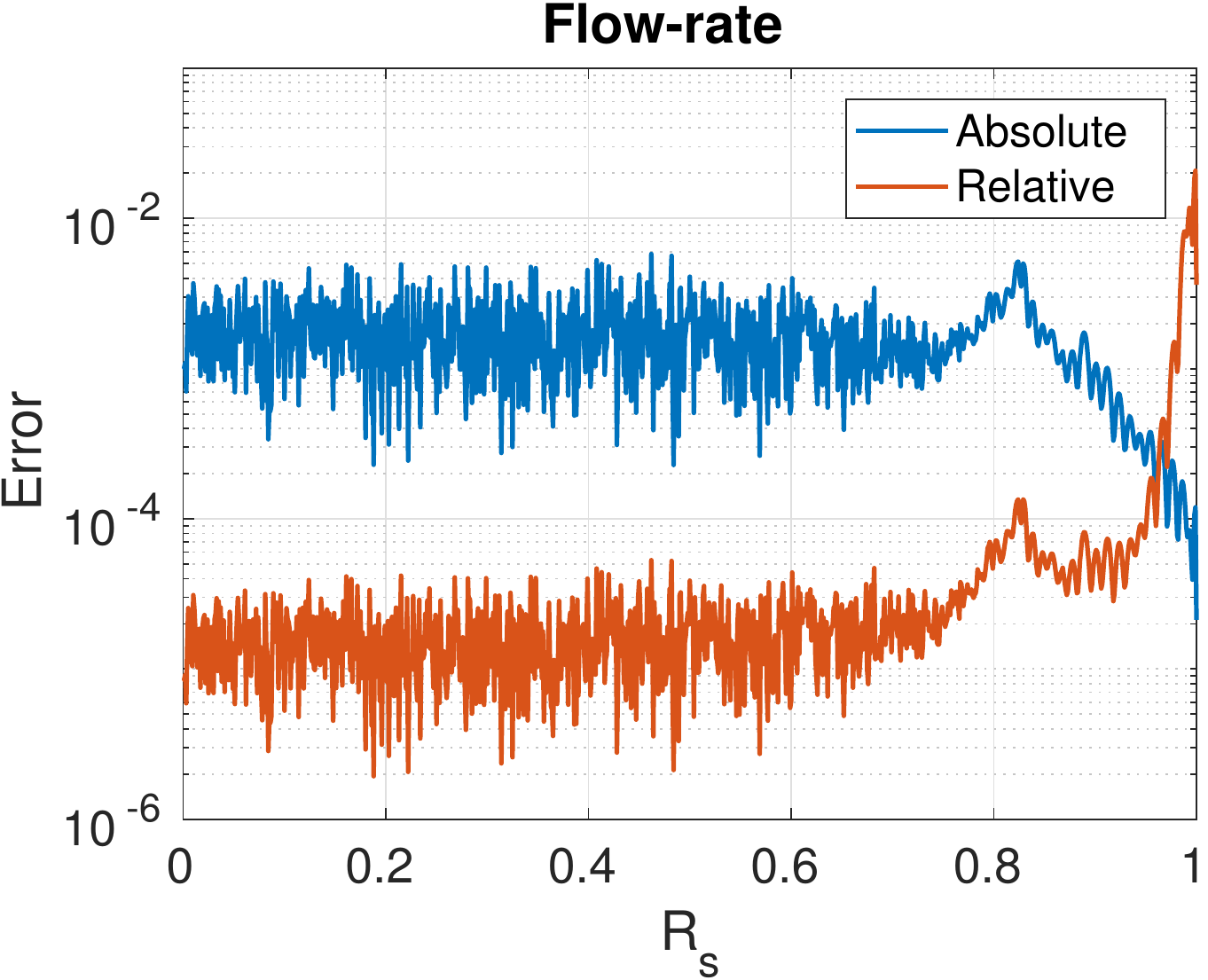}
\end{tabular}
\caption{Pointwise absolute (blue) and relative (red) errors obtained by the kernel models with $N=640$ training points for the 
prediction of the pressure (left) and flow rate (right).}\label{fig:pointwise_errors56}
\end{center}
\end{figure}

We start by describing in detail the results obtained for Vessel 56, i.e., for the functions $f_{56}^P$ and $f_{56}^F$.
Figure \ref{fig:errors56} shows the absolute and relative errors $E^A$ (left), $E^R$ (right) for the two functions. For both the pressure and flow-rate, 
it is clear that an 
increase in the dataset size, hence in the number of full model runs, produces significantly more accurate models. The actual magnitude in the absolute 
errors is different between 
pressure and flow rate, due to a different magnitude of the output quantities. Nevertheless, the relative errors demonstrate that the pressure curves are 
better approximated by the kernel models by about two orders of magnitude, for each dataset size.
In any case, the models exhibit a converging behaviour towards the full model. Moreover, it should be noted that already a relatively 
small dataset of $N=160$ stenosis degrees produces good results in both cases. 
For a better understanding of the error behaviour, we report in Figure \ref{fig:pointwise_errors56} the pointwise absolute and relative errors 
$e_A^{(i)}, e_R^{(i)}$ for pressure (left) and 
flow rate (right) in the case $N=640$. 
Observe that the errors are in all cases very oscillating, since each point in the plots represent one value of 
$e_A^{(i)}$ (or $e_R^{(i)}$) for a different value of $(R_s)_i$, i.e., it is the $2$-norm of the $400$ dimensional vector $f((R_s)_i) - {\hat f}((R_s)_i)$. Thus, the 
small oscillations for a single parameter of the surrogate around the exact solution are amplified into the values depicted in the figures. 

It is worth noticing that in 
both cases the magnitude of the exact quantity is decreasing for $R_s$ close to one, thus the 
relative  error is magnified for high stenosis degrees. This effect is particularly evident for the flow-rate, where the kernel model 
has a better absolute accuracy for 
$R_s\approx1$,  but a significantly worse relative accuracy.

We remark that this causes the worse relative error for the flow rate observed in Figure \ref{fig:errors56}. Indeed, the relative error $E_R$ for a 
given size of the training set is computed as the maximum of the relatives error $e_R^{(i)}$ for all the values $(R_s)_i$ in the test set. Thus, $E_R$ is dominated by 
the 
relative error obtained in the region $R_s\approx 1$, which is large due to the small magnitude of the flow rate computed by the full model. A different error measure, 
e.g. the average relative error, would result in a different error decay in Figure \ref{fig:errors56}.

Similar results have been obtained for the other three vessels. The results are reported in Table \ref{tab:error_table}. It is relevant to notice that 
for these vessels the 
effect of the stenosis is much less visible. Indeed, the relative error for the flow rate is not significantly worse than the one of the pressure, 
since the flow does not completely vanish for increasing stenosis degrees. This effect is more evident for the Vessels $52, 55$, which are not directly 
connected with the stenosis.
\ \\
\begin{table}[h!]
\centering
\begin{footnotesize}
\caption{\label{tab:error_table}  Accuracy of the surrogates for pressure and flow rate for the Vessels $52, 54, 55$. The table shows the absolute ($E_A$) 
and relative ($E_R$) errors obtained with 
datasets of increasing size $N$. }
\begin{tabular}{||c||c|c||c|c||c|c||}
\hline\hline
\multicolumn{7}{||c||}{Pressure}\\
\hline\hline
& \multicolumn{2}{c||}{Vessel $52$} & \multicolumn{2}{c||}{Vessel $54$} & \multicolumn{2}{c||}{Vessel $55$}\\
\cline{2-7}
$N$& $E_A$&$E_R$ & $E_A$&$E_R$ & $E_A$&$E_R$\\
\hline\hline

$10$ & $5.76\cdot 10^{0}$ & $3.01\cdot 10^{-3}$&           $7.35\cdot 10^{0}$ & $3.80\cdot 10^{-3}$           &   $7.41\cdot 10^{0}$ & $3.85\cdot 
10^{-3}$\\
\cline{1-7}
$40$ & $1.15\cdot 10^{0}$ & $6.00\cdot 10^{-4}$                    &$1.50\cdot 10^{0}$ & $7.75\cdot 10^{-4}$             &$1.51\cdot 10^{0}$ & $7.85\cdot 10^{-4}$\\
\cline{1-7}
$160$ & $9.53\cdot 10^{-1}$ & $4.97\cdot 10^{-4}$                &$1.23\cdot 10^{0}$ & $6.35\cdot 10^{-4}$             & $1.25\cdot 10^{0}$ & $6.47\cdot 10^{-4}$\\
\cline{1-7}
$640$ & $2.95\cdot 10^{-1}$ & $1.54\cdot 10^{-4}$                &$2.98\cdot 10^{-1}$ & $1.54\cdot 10^{-4}$                 &$3.01\cdot 10^{-1}$ & $1.56\cdot 10^{-4}$\\
\hline\hline
\multicolumn{7}{||c||}{Flow rate}\\
\hline\hline

$10$ & $6.05\cdot 10^{0}$ & $3.35\cdot 10^{-2}$                & $5.48\cdot 10^{0}$ & $2.51\cdot 10^{-1}$             &$8.33\cdot 10^{-1}$ & $1.20\cdot 
10^{-2}$\\
\cline{1-7}
$40$ & $7.72\cdot 10^{-1}$ & $4.28\cdot 10^{-3}$&                  $5.23\cdot 10^{-1}$ & $2.41\cdot 10^{-2}$                 &$1.84\cdot 10^{-1}$ & $2.64\cdot 10^{-3}$\\
\cline{1-7}
$160$ & $6.35\cdot 10^{-1}$ & $3.52\cdot 10^{-3}$               &$4.08\cdot 10^{-1}$ & $1.88\cdot 10^{-2}$                    &$1.41\cdot 10^{-1}$ & $2.03\cdot 10^{-3}$\\
\cline{1-7}
$640$ & $1.73\cdot 10^{-1}$ & $9.61\cdot 10^{-4}$             &$7.70\cdot 10^{-2}$ & $3.55\cdot 10^{-3}$            &$3.02\cdot 10^{-2}$ & $4.34\cdot 10^{-4}$\\
\hline\hline 
\end{tabular}
\end{footnotesize}
\end{table}

\begin{figure}[ht]
\begin{center}
\begin{tabular}{cc}
\includegraphics[scale=0.38]{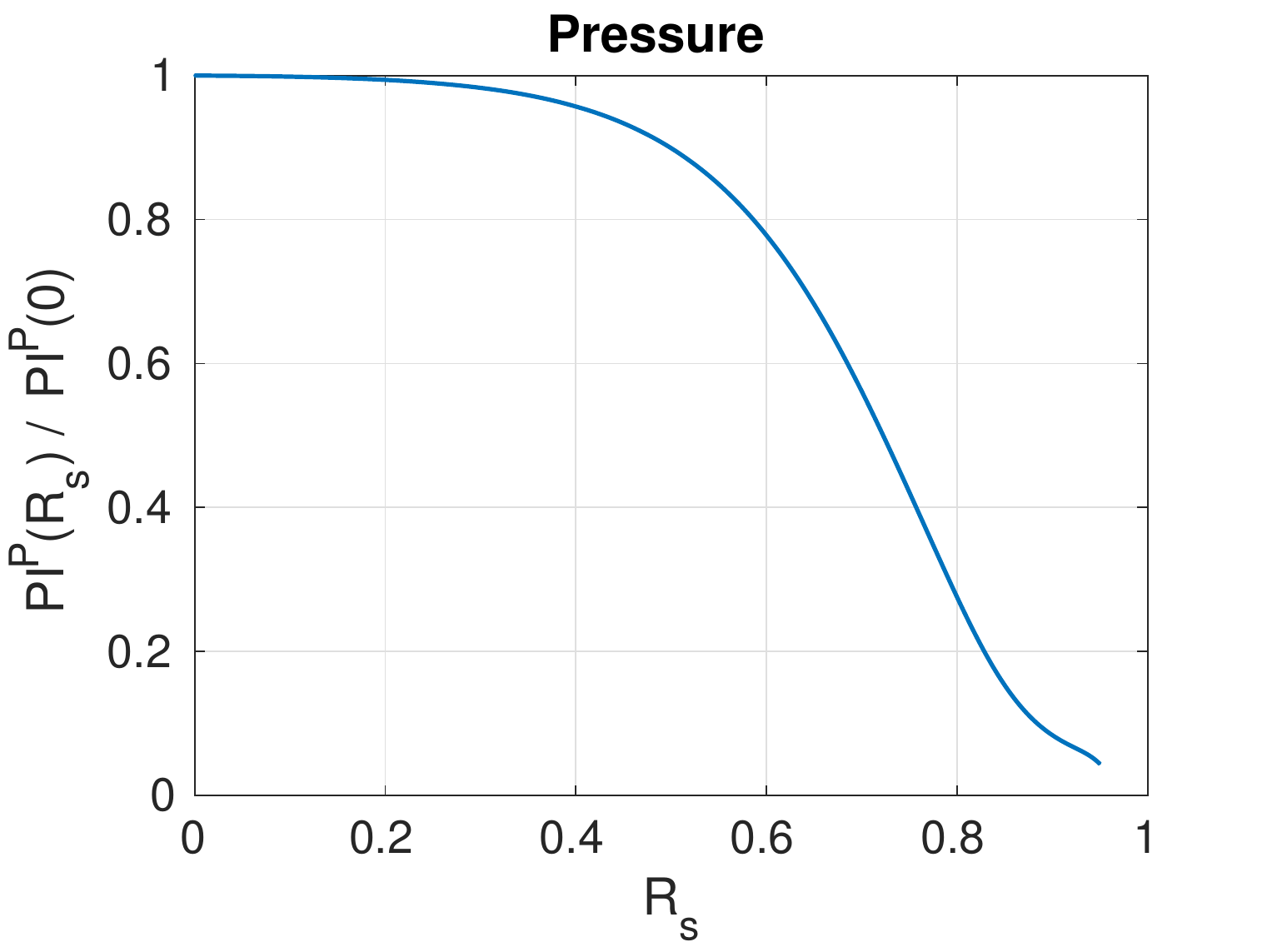}&
\includegraphics[scale=0.38]{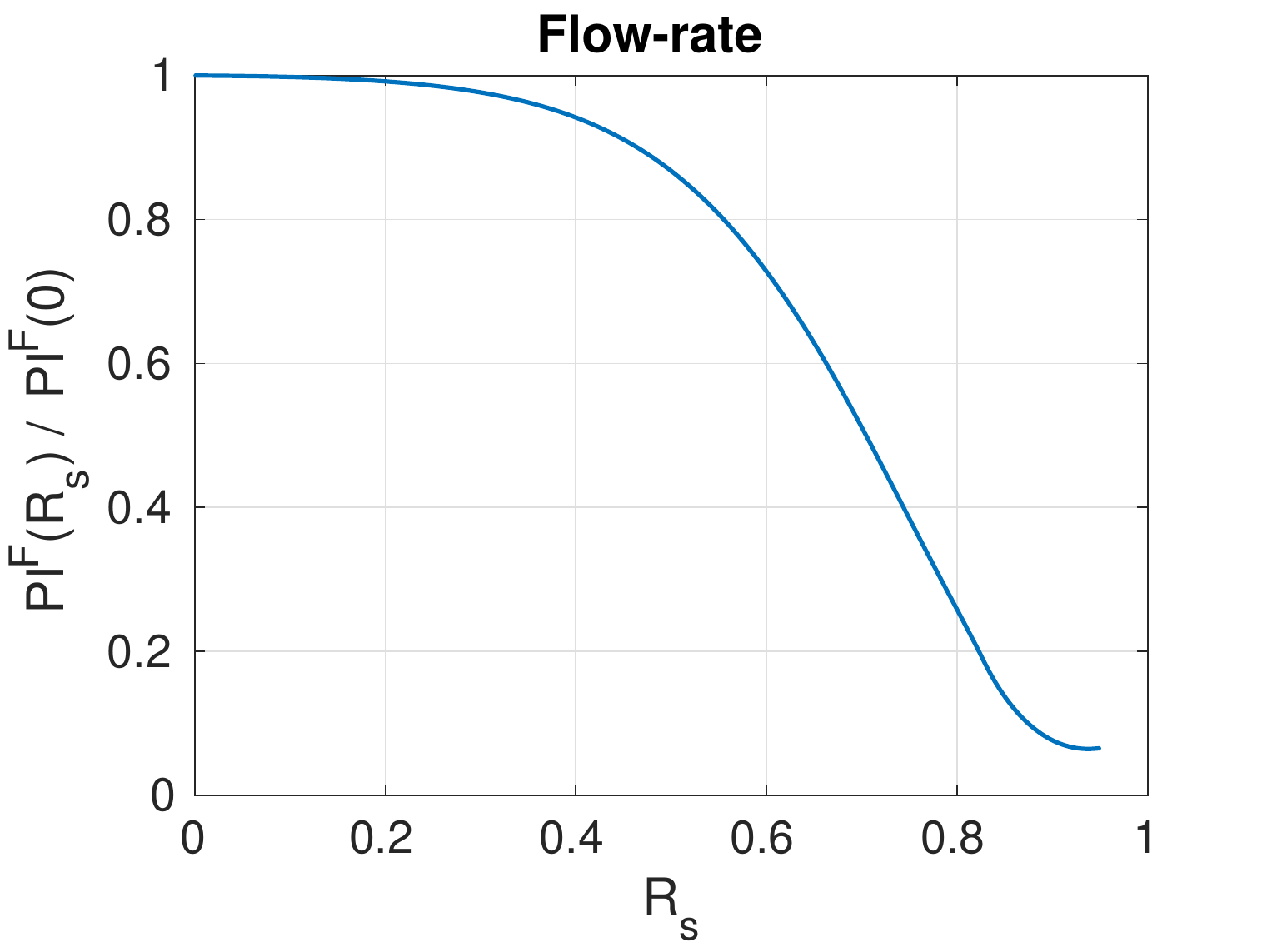}\\
\includegraphics[scale=0.38]{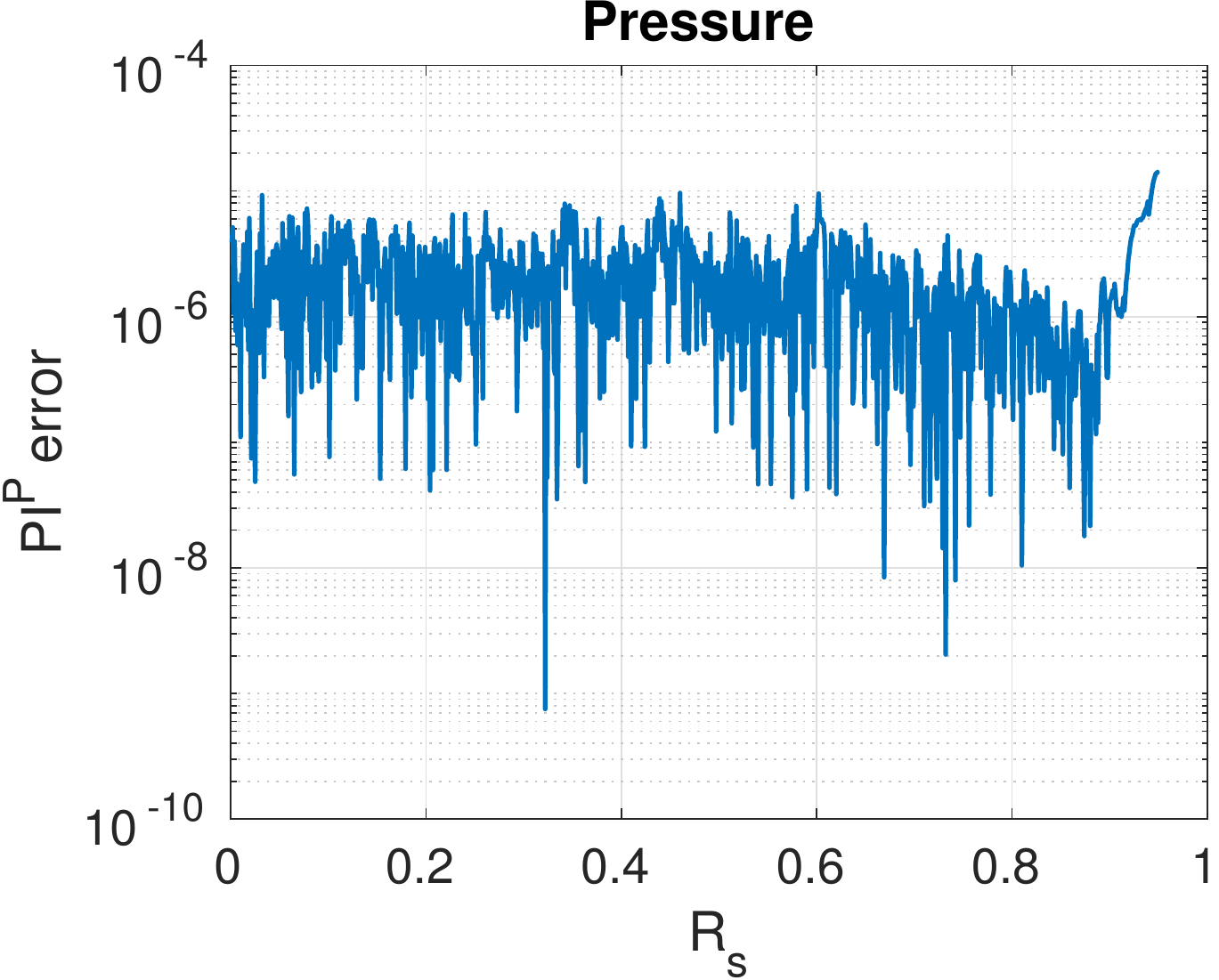}&
\includegraphics[scale=0.38]{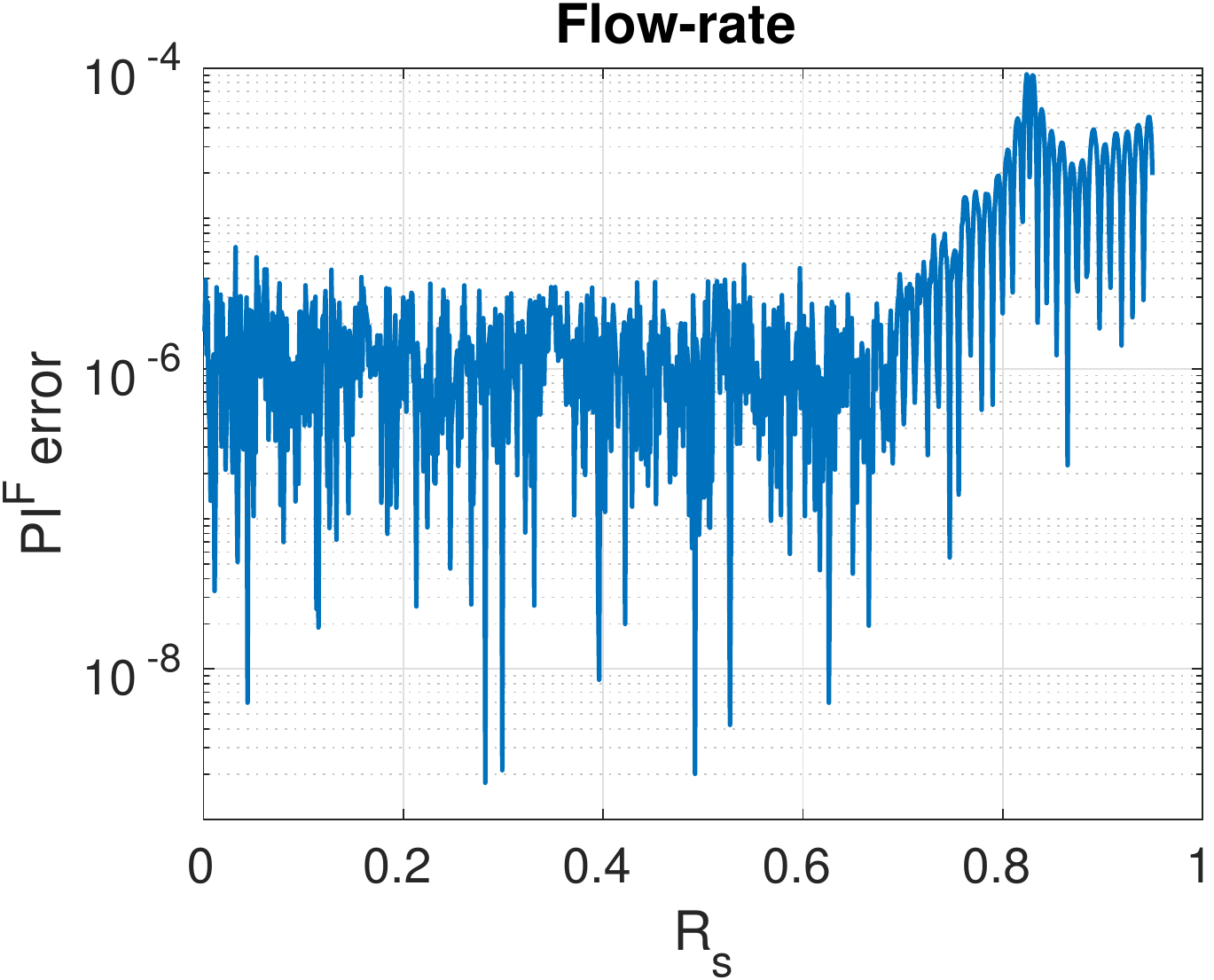}
\end{tabular}
\caption{\label{fig:PI_errors56} Absolute error between the PI computed with the full model and the one computed with the surrogate for the pressure 
(top, left) and flow rate (top, right). The results at the bottom line are obtained by a kernel model with $N=640$.}
\end{center}
\end{figure}

To further measure the accuracy of the surrogate, we compare a relevant blood flow index obtained with the full model and with the surrogate. We consider 
the pulsatility index $PI$, which is a commonly used diagnostic index, and has the advantage of being measurable in a non invasive way (see e.g. 
\cite{stergiopulos1992computer}). For a given stenosis degree $R_s$, in vessel $N_v$ it is computed as
\begin{align*}
PI^P(R_s):=  \frac{\max\left(f_{N_v}^P(R_s)\right) - \min\left(f_{N_v}^P(R_s)\right)}{\mathrm{mean}\left(f_{N_v}^P(R_s)\right)},\;\\
PI^F(R_s):=  \frac{\max\left(f_{N_v}^F(R_s)\right) - \min\left(f_{N_v}^F(R_s)\right)}{\mathrm{mean}\left(f_{N_v}^F(R_s)\right)},
\end{align*}
for the pressure and flow rate, respectively. It measures the difference between the systolic and diastolic pressure (or flow rate) divided by its average value. 
This index can be computed in the same way also for the surrogate using the same formula with the full model replaced with the data-based prediction. Figure 
\ref{fig:PI_errors56} reports the logarithmic error between the exact and the predicted values of $PI^P$ (left) and $PI^F$ (right) in Vessel $56$, for all values of 
$R_s$ in the test 
set, and for the surrogate obtained with $N=640$. The results demonstrate the accuracy of the approximate model also in capturing a physically relevant quantity, and 
this is obtained by data accuracy only, i.e., no constraint is imposed to the surrogate to match the desired values of the $PI$ index of the full model. The simulation
results that are provided in Figure \ref{fig:PI_errors56} show a similar curve progression as in \cite[Figure 3 or Figure 6]{stergiopulos1992computer}. For small stenosis 
degrees, the 
normalised $PI$ indices with respect to $R_s=0$, i.e. the health state, are almost one. As the stenosis degree increases (see Figures \ref{fig:example_pressure_curves} 
and
\ref{fig:example_flow_rate_curves}), it can be observed that in particular the systolic pressures and flow rates are damped significantly. As a consequence, the $PI$ 
indices are decreasing 
monotonously. Up to a stenosis degree of $R_s = 0.5$ the gradient of the curve is rather low, while for a severe stenosis $R_s > 0.5$ the gradient of the $PI$ curve is 
very high.

\subsection{Efficiency of the surrogate models}
It is now of interest to investigate the online efficiency of the surrogates, that is, the time needed to evaluate the models on a new input stenosis. 
In order to understand the timing 
results, we remark that the evaluation of the full model for a given stenosis degree $R_s$ takes around $t_{full} = 200 s$, only for the simulation of the 
time interval [$29s, 
30s]$, i.e., without considering the transient phase. We remark that the present MATLAB implementation can be significantly optimised for speed, and a smaller 
execution time can be expected using a compiled language.

\begin{figure}[ht]
\begin{center}
\begin{tabular}{cc}
\includegraphics[scale=0.4]{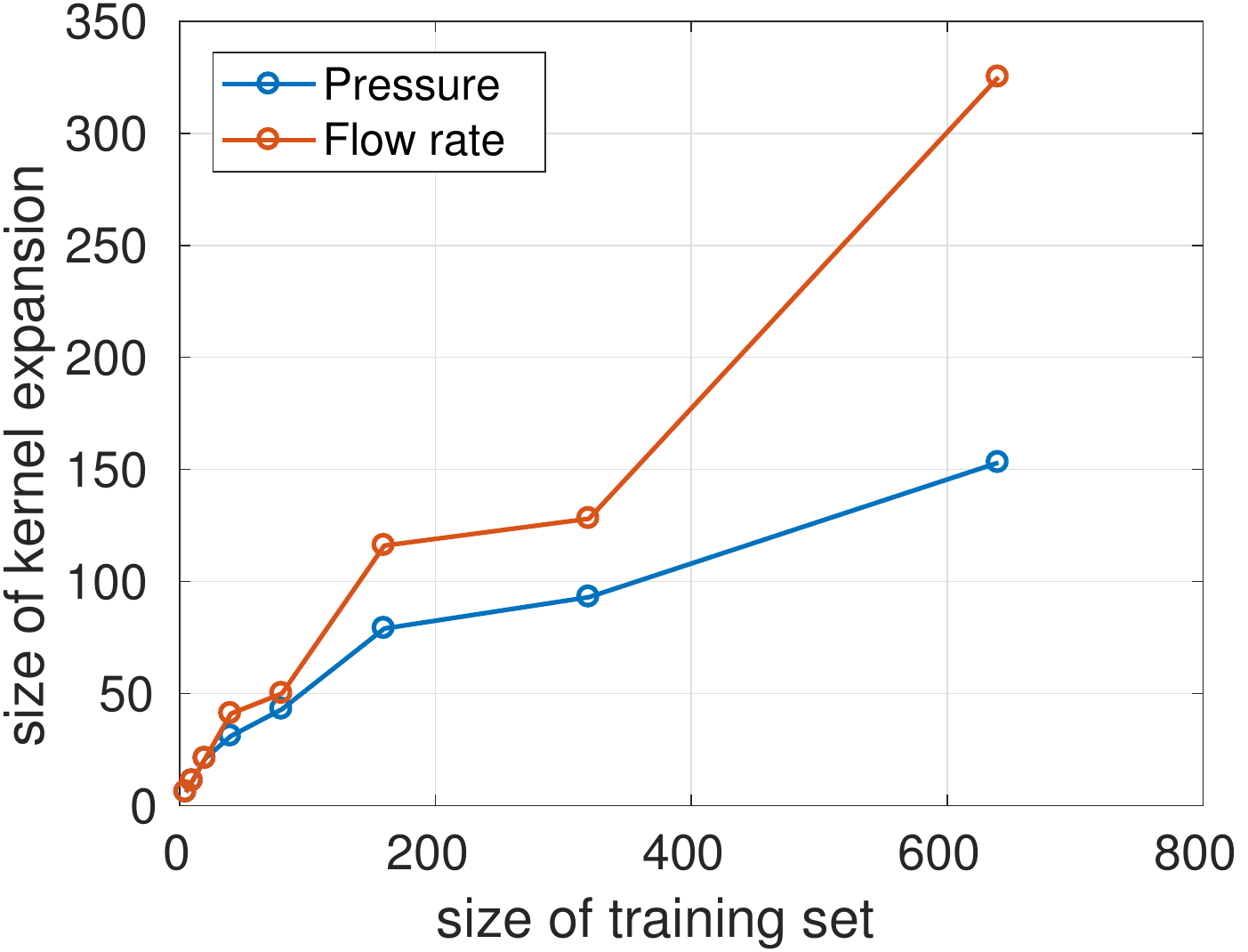}&
\includegraphics[scale=0.4]{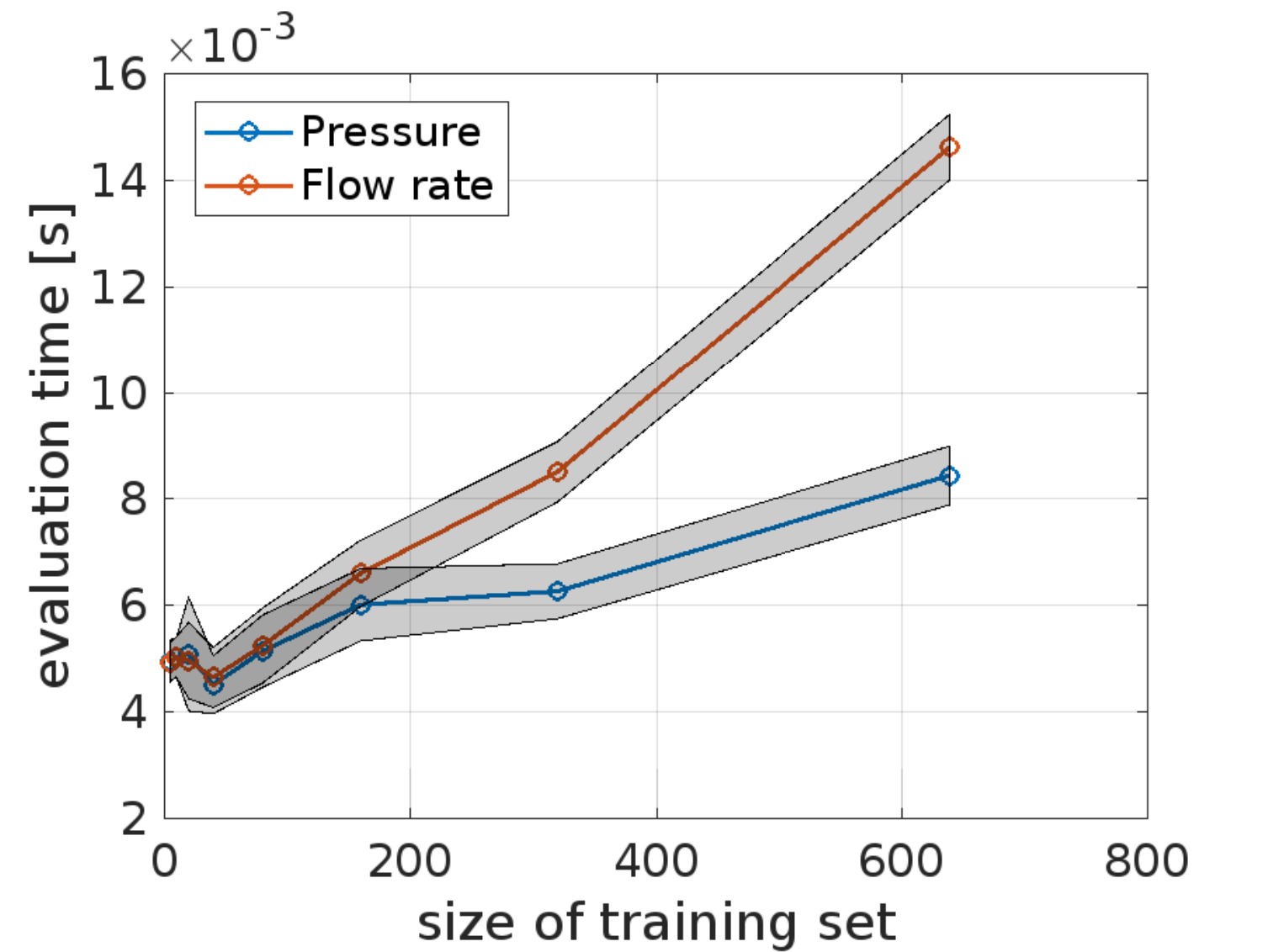}
\end{tabular}
\caption{\label{fig:timing} Number of points selected by the $f$-VKOGA (left) and runtime to evaluate the surrogates on $1000$ inputs (right), for 
Vessel 56 and for datasets of increasing 
sizes.}
\end{center}
\end{figure}

Since the $f$-VKOGA constructs the surrogate selecting only a relevant subset of the full dataset, we look at the actual number of points which is selected for 
the various datasets and output quantities. Figure \ref{fig:timing} (left) reports this number of points, as obtained in the approximation of $f^P_{56}, f^F_{56}$ as 
discussed in the previous section. It is interesting to observe that the number of points increases as the dataset increases, but the number is well below the total 
number of 
points. This means that the surrogates are faster than a non-sparse kernel expansion. Moreover, the flow-rate requires the selection of more points, which confirms that 
this 
output quantity is more difficult to predict. To assess the actual efficiency of the models, in Figure \ref{fig:timing} (right) we report the runtime required to 
evaluate 
each model on the $1000$ test stenosis degrees. The evaluation is repeated $100$ times, and the figure shows the mean and standard deviation over the $100$ experiments.
As expected, the evaluation times are related to the sparsity of the models. 

In all cases the evaluation of the surrogates on $1000$ inputs takes on average less than $1.5\cdot 
10^{-2} s$, i.e., for the largest (hence slowest) model we can estimate an evaluation time per stenosis degree to $t_{surrogate} = 1.5\cdot 10^{-5} 
s$. Compared to $t_{full} = 2 \cdot 10^{2}s$, this still gives a 
speedup factor 
of about $10^{6}$ in the worst case.

\subsection{Solving a state estimation problem by means of kernel methods}
\label{sec:StateEstProb}
In order to demonstrate the use of the surrogate model, we employ it to solve a state estimation problem which would be infeasible by using the full model.
Namely, for a given pressure or flow rate curve in a time interval, we want to predict as accurately as possible the stenosis degree which corresponds to a given curve. 
We assume to collect the measurements into a vector $y\in\mathbb R^q$, which is given by the model output for a fixed, but unknown input stenosis degree $R^{\star}$, 
plus 
some 
additive noise, i.e., 
$$
y := f(R^{\star}) + \eta v,
$$
where $\eta\geq 0$ is a noise level and $v$ is a random vector. In the following experiments, $v$ is drawn from the uniform distribution in $(0, 1)$.
We aim at detecting the value $R^{\star}$ from $y$ only. Doing so, we define a cost function $J(R_s)$ measuring the squared distance between the measurements and the 
model prediction for a given stenosis degree $R_s$, i.e., 
\begin{equation}\label{eq:cost_function}
J(R_s) : = \frac{1}{2\|y_j\|^2_2} \sum_{j=1}^q \left(y_j - \tilde f_j(R_s)\right) ^2,
\end{equation}
and consider the solution $\tilde R^{\star}$ of the optimisation problem
\begin{equation*}
\tilde R^{\star} = \min_{R_s\in[0,1]} J(R_s). 
\end{equation*}
Observe that, when the noise term $\eta v$ is vanishing and the model prediction is exact, the unique minimiser is the exact solution, i.e., 
$\tilde R^{\star} = R^{\star}$.

In principle, it would be possible to formulate the cost function \eqref{eq:cost_function} also in terms of the full model $f$. 
Nevertheless this is infeasible in practice, since 
multiple evaluations of $f$ are required to compute a minimiser. The use of the cheap surrogate model, instead, allows a real-time estimation of $R^{\star}$.
Moreover, since the kernel is differentiable the cost function is also differentiable, thus the use of 
gradient-based methods is possible. In particular, we have 
\begin{equation*}
\frac{d}{dR_s}J(R_s)  = \frac{1}{\|y_j\|^2_2}\sum_{j=1}^q \left(y_j - \tilde f_j(R_s)\right)\frac{d}{dR_s}f_j(R_s),
\end{equation*}
with 
\begin{equation*}
\frac{d}{dR_s} f_j(R_s) = \dd{R_s}{\sum_{i=1} ^N \alpha_{ij} K(R_s, R_i)} = \sum_{i=1} ^N \alpha_{ij} \frac{d}{dR_s} K \left(R_s, R_i \right)
\end{equation*}
and the $R_s$-derivative of the kernel can be explicitly computed, since the kernel itself is known in closed form. 
In other terms, both the cost function $J$ and its derivative can be computed efficiently by means of the surrogate, and they both only involve 
the evaluation of matrix-vector products. 

\begin{figure}[ht]
\begin{center}
\begin{tabular}{cc}
\includegraphics[scale=0.4]{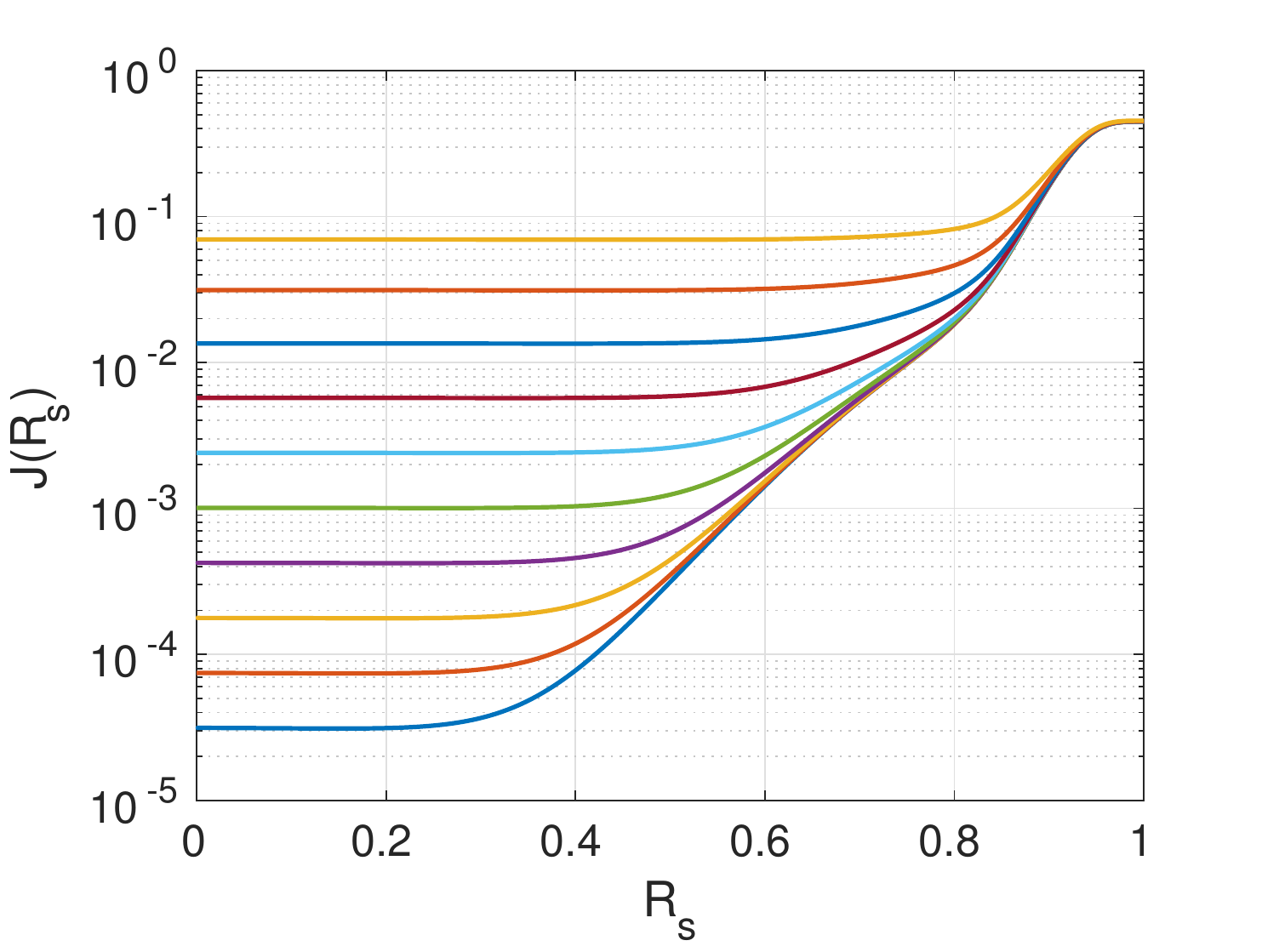}&
\includegraphics[scale=0.4]{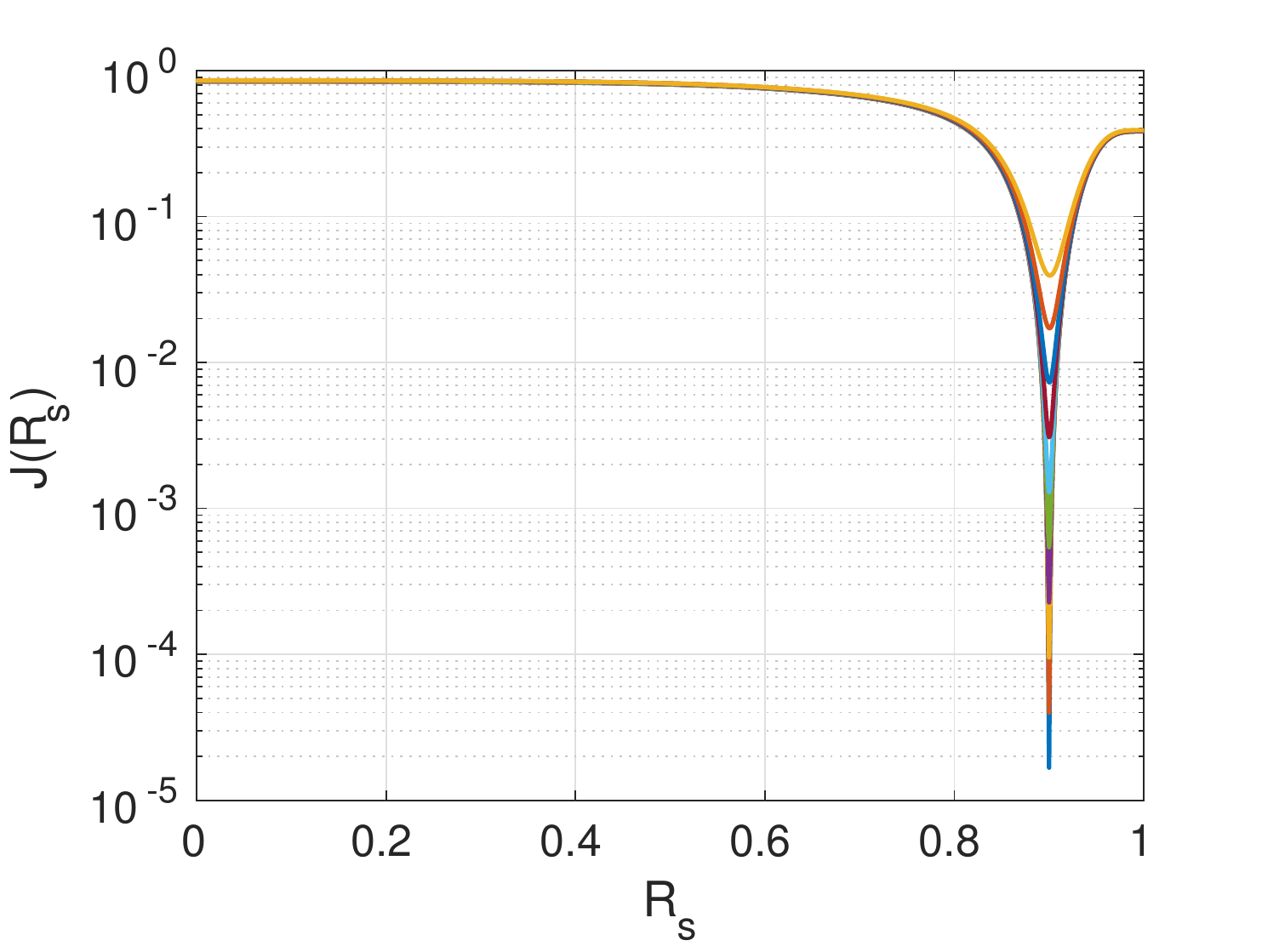}\\
\includegraphics[scale=0.4]{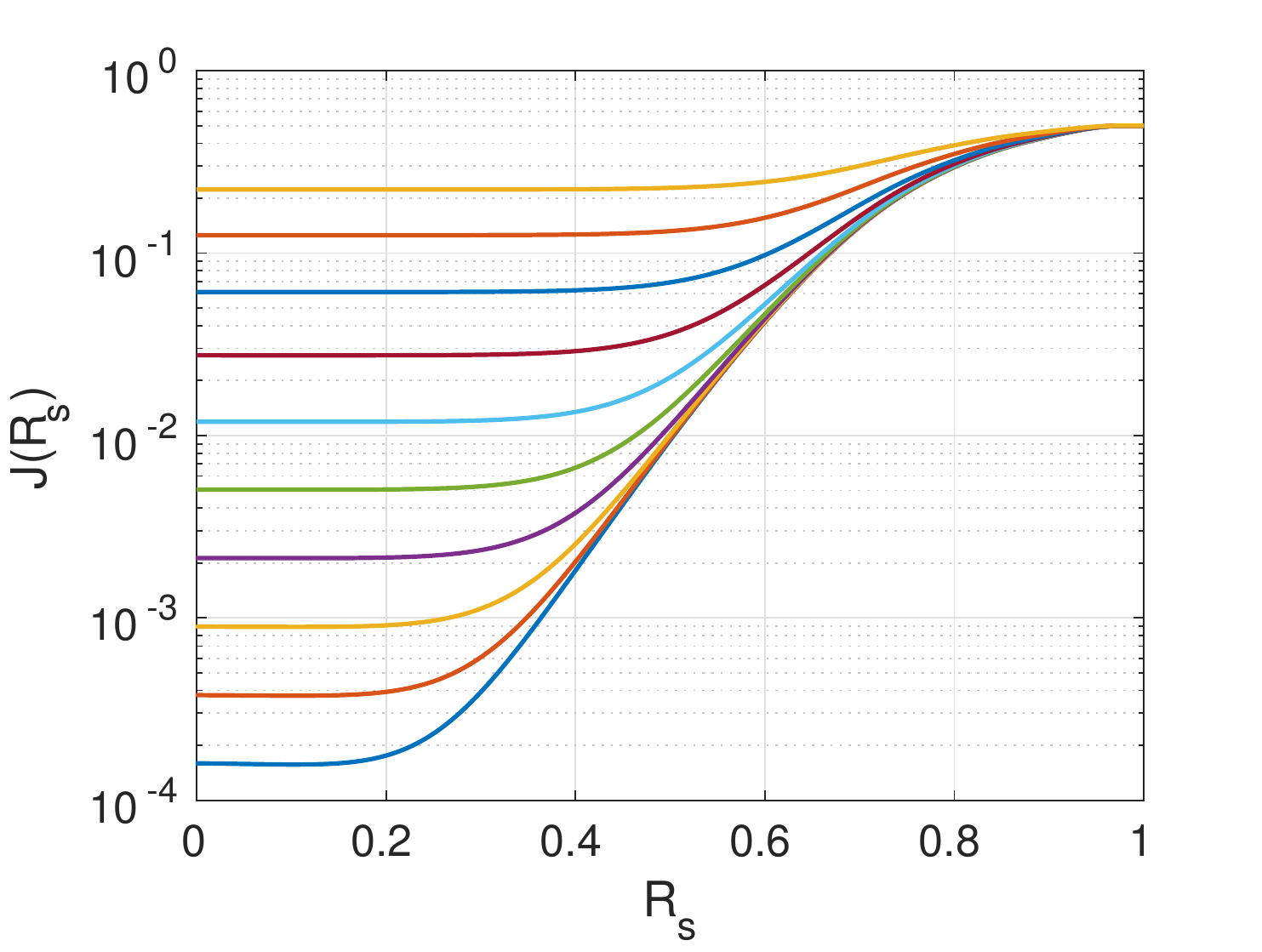}&
\includegraphics[scale=0.4]{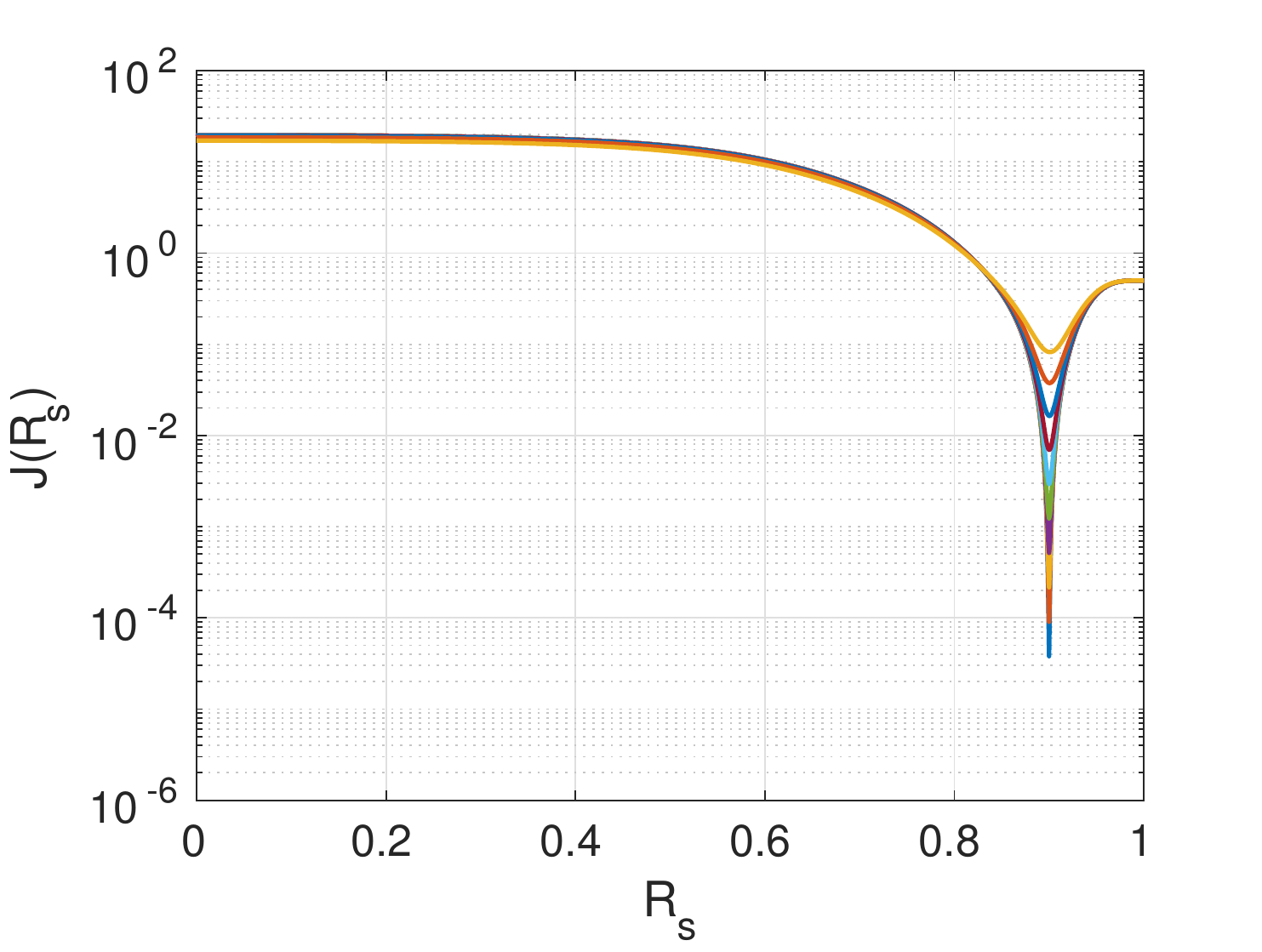}
\end{tabular}
\caption{\label{fig:cost_functions} Cost functions $J$ for $R^{\star} = 0.1$ (left) and $R^{\star} = 0.9$ (right) for the pressure curve (top line) 
and the flow rate (bottom line), for Vessel 56 and using the surrogate based on the dataset of $N=160$ full model runs. The various curves correspond to 
increasing noise levels in the target value $y$.}
\end{center}
\end{figure}

We proceed to some experiments for Vessel 56 and the surrogate models obtained with the dataset of $N=160$ full model runs. 
In Figure \ref{fig:cost_functions}, we plot the cost function $J$ for $R_s\in[0,1]$ and $R^{\star} = 0.1$ (left) and $R^{\star}= 0.9$ (right) 
for the pressure (top) and flow rate 
(bottom). The results are reported for increasing noise levels $\eta$, obtained as $10$ logarithmically spaced values in $\left[0.01,0.5 \right]$. It is 
clear that the surrogate provides a 
reliable prediction when the stenosis degree is large, also in the presence of noise, since the cost function has a unique minimiser. Instead, for small target 
stenosis degrees the cost function is flat, so we should expect a less accurate prediction. 

\begin{figure}[ht]
\begin{center}
\begin{tabular}{cc}
\includegraphics[scale=0.41]{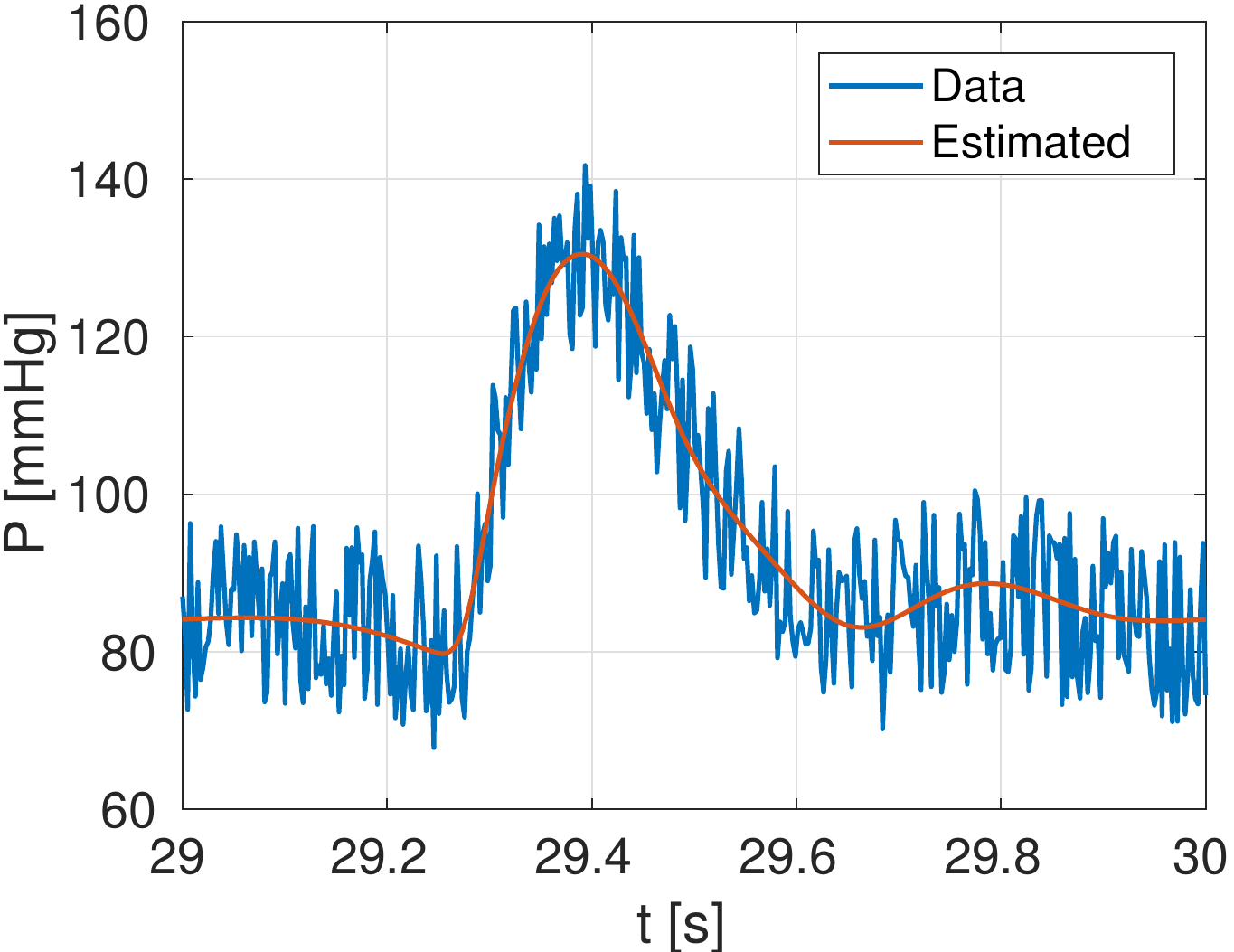}&
\includegraphics[scale=0.41]{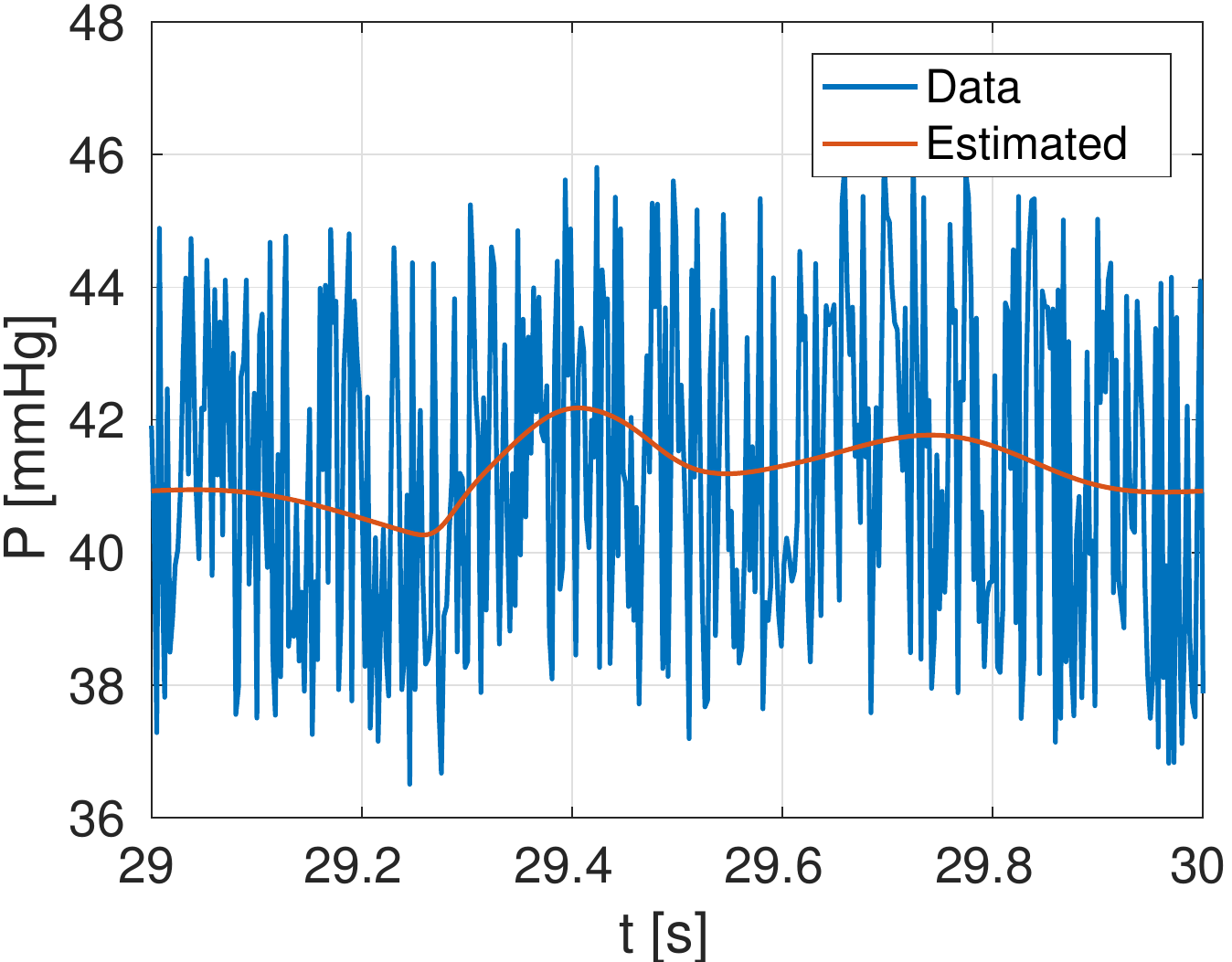}\\
\includegraphics[scale=0.41]{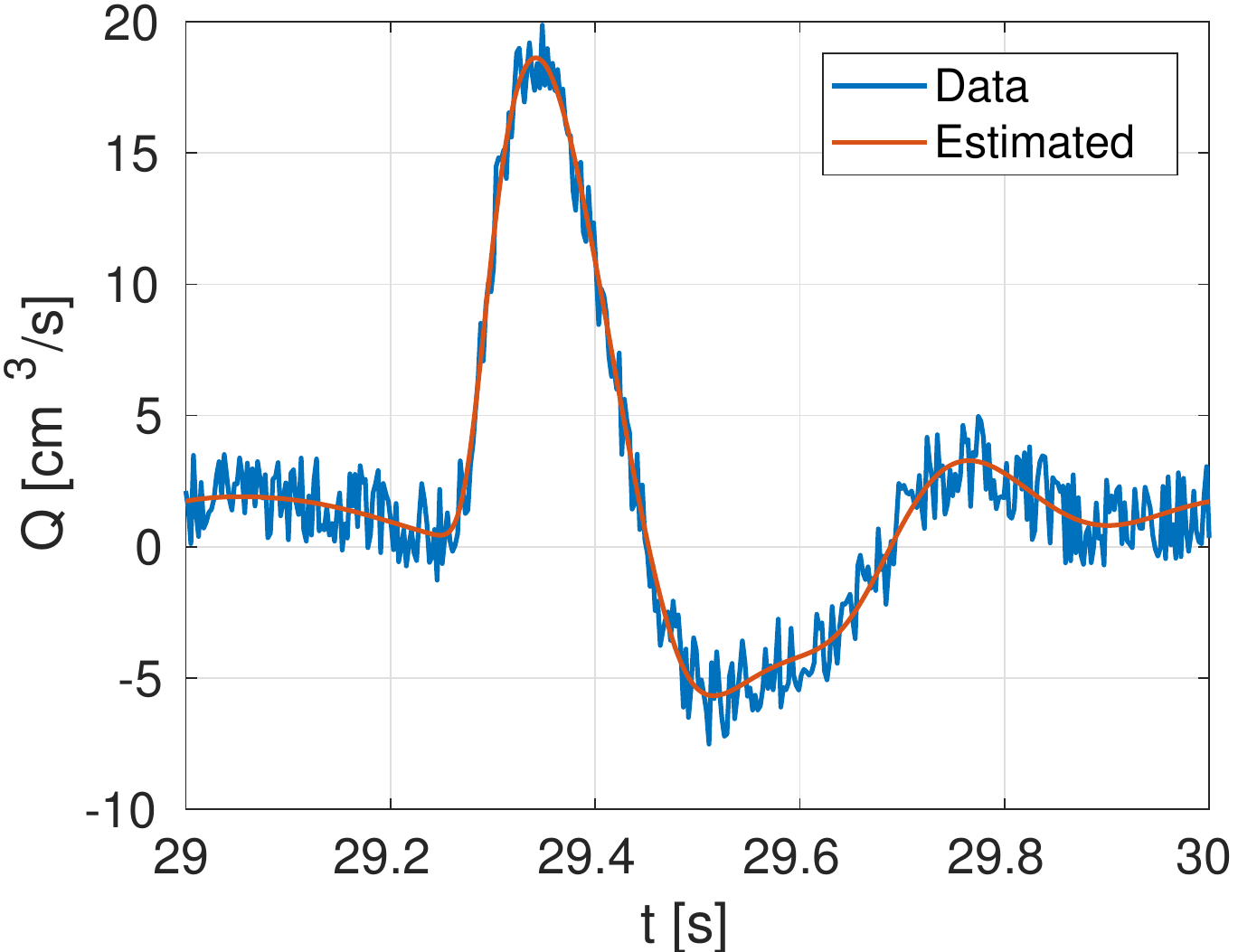}&
\includegraphics[scale=0.41]{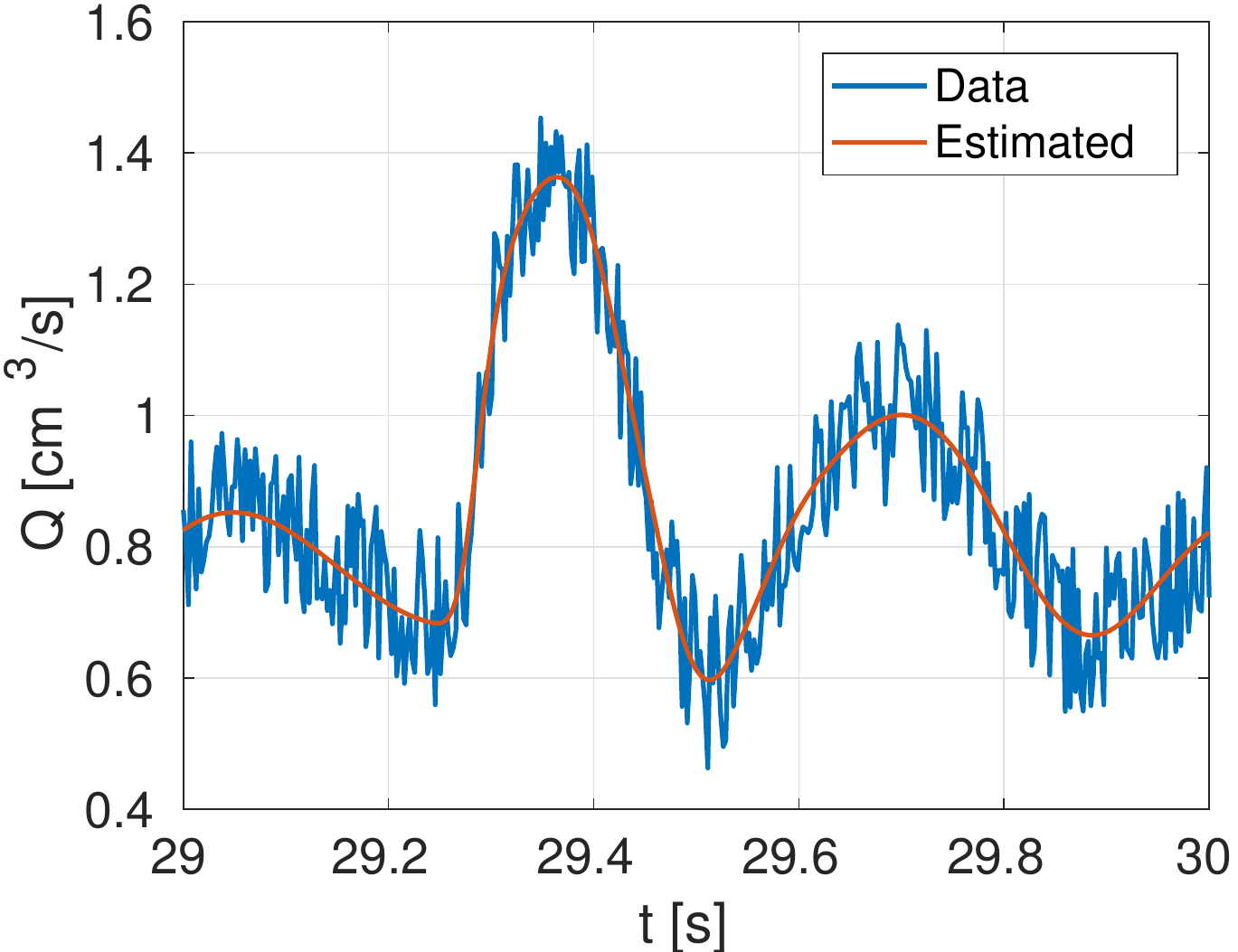}
\end{tabular}
\caption{\label{fig:estimations} Noisy data and estimated curves for $R^{\star} = 0.1$ (left) and $R^{\star} = 0.9$ (right) for the pressure 
curve (top line) and the flow rate (bottom line) (Vessel 56). The surrogate model is based on the dataset of $N=160$ full model runs.}
\end{center}
\end{figure}

The values of $\tilde R^{\star}$  can then be computed with any constrained optimisation solver, and we use here the MATLAB built-in \texttt{fmincon}, which uses an 
active set search procedure. The noisy input $y$ for $\eta = 0.1$, as well as the resulting estimated curves, are depicted in Figure \ref{fig:estimations}, for 
both $R^{\star}=0.1, 
0.9$ and both pressure and flow rate. The estimated values $\tilde R^{\star}$ are reported in Table \ref{tab:optim_error}.

\ \\
At a first look, the results could seem somehow surprising, since the estimation is much better for large stenosis degrees, i.e., in the cases where the surrogates are less 
accurate. Nevertheless, the exact values of the pressure and the flow rate indeed are less variable for small $R_s$, so for the estimator it is more difficult to discriminate between 
different values.

\begin{table}[h!]
\centering
\caption{\label{tab:optim_error} Results of the parameter identification problem for the pressure and flow rate curves with 
$R^{\star}=0.1$ and $R^{\star} = 0.9$. The table reports the estimated values $\tilde 
R ^{\star}$ and the errors with respect to the exact state.}
\begin{tabular}{||c||c|c||c|c||}
\hline\hline
& \multicolumn{2}{c||}{Pressure} &\multicolumn{2}{c||}{Flow rate}\\
\hline
$R^{\star}$& $\tilde R ^{\star}$ & Error & $\tilde R ^{\star}$ & Error\\
\hline
$0.1$& $3.11 \cdot 10 ^{-1}$&$2.11 \cdot 10 ^{-1}$   &$1.37 \cdot 10 ^{-1}$&$3.71 \cdot 10 ^{-2}$\\
\hline
$0.9$& $9.00 \cdot 10 ^{-1}$&$1.67\cdot 10^{-4}$  &$9.00 \cdot 10 ^{-1}$&$1.83\cdot 10^{-4}$\\
\hline\hline
\end{tabular}
\end{table}

\section{Conclusion and further work}
\label{sec:Concl}

In this paper, we have simulated blood flow in the $55$ main arteries of the systemic circulation. For this purpose 1-D blood flow
models have been considered. At the outlets of the main arteries $0$-D \\
lumped parameter models have been coupled with the corresponding
$1$-D models to include the Windkessel effect of the omitted vessels, while at the inlet of the aorta (Vessel 1) a lumped parameter model for the
left ventricle has been coupled with the $1$-D model for the aorta. Furthermore, the impact of a stenosis in a tibial artery has been simulated
by an ODE depending on the degree of the stenosis. In order to be able to obtain insight into the flow behaviour in the vicinity of the stenosis
for an arbitrary degree of stenosis without starting the simulation for every degree of stenosis again, the output data are provided by a surrogate
model based on kernel methods. It is has been demonstrated that the error between the surrogate kernel model and the exact 
simulation result is decreasing as more and more training data are included into the surrogate model. In addition to that the efficiency of the
kernel method has been investigated. The surrogate kernel model has been used to solve a parameter and state estimation problem: For a given pressure
or flow rate curve in the vicinity of the stenosis the corresponding degree of stenosis is estimated, and the kernel-based surrogate yields also a prediction 
of the state. This is a step towards real-time estimation and decision in patient-specific treatments.

Future work may be concerned with simulating the whole circulation by means of a closed loop model. This means that besides the 
left ventricle also the remaining chambers of the heart as well as the pulmonary circulation and the venous part of the systemic
circulation have to be modelled \cite{liang2009multi}. A further aspect that could be investigated would be to 
include besides the stenosis degree more parameters of interest. Two possible parameters would be the peripheral resistances of Vessel $55$ and $56$. In Subsection 
\ref{sec:PerCirc}, we denoted them by $R_{2,55}$ and $R_{2,56}$. Varying these
parameters, one could simulate the effect of vasodilation, i.e. the enlargment of arterioles due to a reduced blood supply of tissue.
In this context it is of great interest, how these resistances have to be adapted
such that for a given degree of stenosis, a maximal blood flow rate distal to the stenosis can be restored.

Moreover, the cost function for the state estimation problem in Subsection \ref{sec:StateEstProb} could be improved such that the estimates for low stenosis
degrees are more accurate.

In more general terms, other aspects of model- and data-driven surrogate modeling could be investigated. For instance, the same kernel based technique has 
been recently applied to Uncertainty Quantification \cite{koeppel2017}. In the present setting, this could lead to the fast assessment of the impact of uncertainty on 
the model output, e.g. in a setting where the stenosis degree or other possible input quantities are not exactly measured.

\section*{Acknowledgements}
This work was supported by the Cluster of Excellence in Simulation Technology (EXC 310/2) at the University of Stuttgart. 

\bibliographystyle{plain}

\nocite{*}

\end{document}